\documentclass[12pt,english]{article}
\usepackage[T1]{fontenc}
\usepackage[latin9]{inputenc}
\usepackage[paperwidth=8.5in,paperheight=11in]{geometry}
\usepackage{longtable}
\usepackage{amstext}
\usepackage{amssymb}
\usepackage{setspace}
\makeatletter

\providecommand{\tabularnewline}{\\}

\usepackage{enumitem}       

\usepackage{aStyle}

\makeatother

\usepackage{babel}
\begin{document}
\begin{doublespace}
\begin{center}
\textbf{\Large{}How to Solve Diophantine Equations}\vspace{-1.3cm}
\par\end{center}
\end{doublespace}

\begin{center}
Taha Sochi (Contact: ResearchGate)\vspace{-0.4cm}
\par\end{center}

\begin{center}
London, United Kingdom
\par\end{center}

\textbf{Abstract}: We present in this article a general approach (in
the form of recommendations and guidelines) for tackling Diophantine
equation problems (whether single equations or systems of simultaneous
equations). The article should be useful in particular to young ``mathematicians''
dealing mostly with Diophantine equations at elementary level of number
theory (noting that familiarity with elementary number theory is generally
required).

\begin{flushleft}
\footnote{All symbols and abbreviations in this paper are defined in $\S$ \hyperref[Nomenclature]{Nomenclature}.}
\par\end{flushleft}

\clearpage{}

\tableofcontents{}

\clearpage{}

\section{\label{secIntroduction}Introduction}

The subject of Diophantine equations is one of the oldest and richest
mathematical branches in number theory (which is one of the oldest
and richest branches, if not the oldest and richest, of all mathematics).
The roots of this subject can be traced back to the ancient Babylonians
who dealt with Pythagorean triple problems (which are essentially
Diophantine equations).\footnote{In fact, the credit for the ``Pythagorean triples'' should be attributed
in the first place to the Babylonians (and possibly other civilizations
of ancient Mesopotamia) who recorded and used these triples in their
calculations (mainly for the purpose of surveying farm lands) more
than one millennium before Pythagoras and the Greek Pythagorean school.
There is also evidence (or indication) that the ancient Egyptians
used these triples in their calculations for similar reasons and purposes
as the Babylonians (and possibly even in the engineering of their
marvelous constructions like the great Pyramids). So, it should be
more fair to call these triples something like ``Babylonian triples''
or ``Babylonian-Pythagorean triples''. For more details about this
historical issue, the reader should refer to the literature of the
history of mathematics with special attention to the literature about
the Babylonian clay tablets related to the mathematics of ancient
Mesopotamia (e.g. Plimpton 322 and Si.427 tablets) which are investigated
by a number of archaeologists, mathematicians and historians of science
(such as Dr. Daniel Mansfield of the University of New South Wales;
see for instance \cite{Mansfield2021}).}

The label of these equations as ``Diophantine equations'' is another
indication to their ancient roots since it originates from the name
of the ancient Greek-Egyptian mathematician (i.e. Diophantus of Alexandria)
who lived in the third century and who documented in his books problems
which are modeled by this type of equations. In fact, being a Greek-Egyptian
mathematician should also suggest an oriental origin of these types
of equations and the mathematical branch that is based on them. Apparently,
ancient Hindu mathematicians have also dealt with this type of equations
and the methods of their solution as early as the fifth or sixth century
BC.

Despite their apparent simplicity, Diophantine equations are usually
difficult (if not impossible) to solve. In fact, Diophantine problems
(represented mostly in solving Diophantine equations and systems of
such equations) are generally more difficult to tackle and solve than
their corresponding ordinary versions. This is because the demand
for the solutions (and answers in general) to be integers imposes
extra requirements and conditions and hence it usually complicates
the process and methods of solution. In this context, it is worth
noting that one of the most famous open problems in mathematics (which
waited several centuries of mathematical development and required
huge efforts from many prominent mathematicians before it was finally
solved in the mid 1990's) is Fermat's last theorem which is a Diophantine
equation problem. Also, there are still many unsolved open problems
about the subject of Diophantine equations or related to it.

In this article we propose a general strategy for tackling and solving
Diophantine equations. This strategy is presented in the form of a
list of recommendations and guidelines that can (and should) be used
as a reference in tackling Diophantine problems. In fact, most of
these recommendations and guidelines are based on (and extracted from)
our personal experience in solving (as well as in reading) Diophantine
equation problems and the methods (or techniques or tricks or $\ldots$)
which are commonly used in their solution (see for example \cite{SochiBookNTV1,SochiBookNTV2}).
The structure of this paper is simply based on the aforementioned
list where each one of the following sections is essentially based
on one of the recommendations and guidelines (with consideration of
their order in part, i.e. when order is relevant or required).\footnote{However, we should note that the order of the following sections (and
hence the recommendations they present) does not necessarily reflect
the required order of these recommendations in practical situations
noting that the order generally depends on the type of the Diophantine
problem and its characteristic features. Anyway, commonsense and wise
judgment should always be the resort that determines the order (as
well as almost everything else).}

We should finally note (before we go through this investigation) that
a Diophantine equation (as well as a system of equations) may have
no solution, or a single solution, or multiple solutions (whether
finitely many or infinitely many). So, ``solving'' (or ``finding
the solution'') of a Diophantine equation (or system) should mean
\textit{proving} (by an irrefutable logical/mathematical argument)
that there is no solution (i.e. when there is no solution) or \textit{finding}
all the solutions (either explicitly or through a sort of closed form
formula or formulae) with an incontestable argument that there are
no other solutions. So, a Diophantine problem (whether equation or
system) is not solved, for instance, by finding a number of solutions
(e.g. through inspection or through computational search) even if
we know for sure that there are no other solutions. We should also
note that a solution of a system of Diophantine equations is a solution
that satisfies all the equations in the system simultaneously (i.e.
the set of solutions of a system of Diophantine equations is the intersection
of the sets of solutions of its individual equations). So, these criteria
and considerations about ``solving'' Diophantine equations and systems
should apply throughout this paper.

We should also note that there are two main methods for solving systems
of Diophantine equations in number theory. The first is based on using
the traditional methods of solving systems of multivariate equations
(as investigated in algebra and linear algebra for instance) such
as by substitution or comparison or use of the techniques of matrices,
and the second is by solving the individual equations separately (either
by the general methods of algebra or by the special methods and techniques
of number theory) and selecting the solutions that satisfy the system
as a whole (i.e. by accepting only the solutions which are common
to all the equations).\footnote{In this context, we should pay special attention to the method of
having (or obtaining) the solution of one equation and test it on
the other equations in the system. In brief, if the system of Diophantine
equations contains an equation whose solution is known (or can be
obtained easily or more easily) then the best approach for solving
the system is to test the solutions of that equation on the other
equations where only the common solutions (if any) to all equations
are accepted. This is particularly true when some of the equations
in the (non-linear) system are linear. This approach usually saves
considerable amounts of time and effort in trying to solve the system
by other methods. In fact, in some cases this can be the only viable
method for solving the system.} In this regard we remind the reader of what we already said, that
is: the set of solutions of a system of equations is the intersection
of the sets of solutions of its individual equations. As a result,
a system of equations is solvable only if its individual equations
are solvable, although the converse is not true in general. Accordingly,
a system of equations has no solution if some of its equations have
no solution, but a system may not have a solution even though all
its individual equations have solutions (i.e. when the intersection
of these solutions is the empty set).

\clearpage{}

\section{Initial Sensibility Checks}

The first recommendation is to conduct initial (and basic) sensibility
checks to assess the sensibility of the equation quickly (by inspecting
its general characteristics) to see if it is possible to have a solution
or not. These initial checks may also reveal the obvious solutions
of the equation easily without effort or use of any complicated treatment.
In the following subsections we outline some of the most common initial
sensibility checks (where we dedicate the last subsection to our ``final
thought'' about this issue).

\subsection{\label{secParityChecks}Parity Checks}

Parity checks should be regarded as the first item in the list of
sensibility checks. This is due to its simplicity and intuitivity.
For example, if we are asked to find the general solution of the Diophantine
equation:
\[
x^{4}+4y^{3}-7x^{2}-12y+7=0\hspace{2cm}(x,y\in\mathbb{Z})
\]
then before we try to solve this equation by using the familiar rules
and traditional methods of solving polynomial Diophantine equations
(in two variables) we should simply check the parity of this polynomial,
and hence we can easily conclude (by checking the parity) that this
equation has no solution because the polynomial is always odd and
hence it cannot be equal to 0 which is even.

Similarly, the Diophantine equation:
\[
x^{9}-x=y^{5}-13y^{2}+1\hspace{2cm}(x,y\in\mathbb{Z})
\]
can be ``solved'' with no effort by noting that the left hand side
is always even while the right hand side is always odd and hence this
equation has no solution in integers (as it is supposed to be a Diophantine
equation).

Another example is the exponential equation:
\[
17^{x}-13^{y}=19^{z}\hspace{2cm}(x,y,z\in\mathbb{N}^{0})
\]
which obviously has no solution because the LHS is even while the
RHS is odd.

Parity checks can also reduce the possibilities that to be considered
(or the domain of the problem). For example, the equation:
\[
18^{x}+16^{y}=19^{z}\hspace{2cm}(x,y,z\in\mathbb{N}^{0})
\]
can have a solution (in principle) but because of parity considerations
any potential solution must have either $x=0$ (and $y\neq0$) or
$y=0$ (and $x\neq0$). So, we have only these possibilities to consider
which by simple inspection should lead to the only solution, i.e.
$(x,y,z)=(1,0,1)$.

\subsection{\label{subPrimalityComposityChecks}Primality and Composity Checks}

For example, if we conduct an initial primality check on the following
Diophantine equation:
\[
5x^{2}+125y^{3}=4973\hspace{2cm}(x,y\in\mathbb{Z})
\]
then it should be fairly obvious that this equation has no solution
because 4973 is prime while $5x^{2}+125y^{3}=5(x^{2}+25y^{3})$ and
hence they cannot be equal considering their prime factorization.\footnote{Alternatively, we may say: because 4973 is prime while $5x^{2}+125y^{3}=5(x^{2}+25y^{3})$
which is either composite or equal to 5 and hence they cannot be equal.}

Similarly, the equation:
\[
6x^{3}-19x^{2}y+19xy^{2}-6y^{3}=p\hspace{2cm}(x,y\in\mathbb{Z}\textrm{ and }p\in\mathbb{P})
\]
 has no solution because the left hand side cannot be prime (considering
its factorization) and hence this Diophantine equation has no solution.\footnote{The LHS can be factorized as $(2x-3y)(3x-2y)(x-y)$. Now, a prime
number $p$ can be factorized as a product of three integer factors
only in 3 different ways, i.e.
\[
p=(1)(-1)(-p)=(1)(1)(p)=(-1)(-1)(p)
\]
It can be shown that none of these ways are applicable in this case
and hence the LHS cannot represent a prime number.}

\subsection{Sign and Magnitude Checks}

For example, the Diophantine equation:
\[
x^{2}+y^{2}+1=0\hspace{2cm}(x,y\in\mathbb{Z})
\]
has obviously no solution because $x^{2}+y^{2}$ cannot be negative
and hence when it is added to 1 the result cannot be zero.

Similarly, the equation:
\[
x^{4}+y^{2}+z^{6}=0\hspace{2cm}(x,y,z\in\mathbb{Z})
\]
has the obvious (and only) trivial solution (i.e. $x=y=z=0$) because
any sum of even natural powers of integers must be a positive natural
number unless all the integers are 0.

Also, it is fairly obvious that the equation:
\[
\frac{1}{x}+\frac{1}{y}+\frac{1}{z}=5\hspace{2cm}(x,y,z\in\mathbb{Z}\textrm{ and }xyz\neq0)
\]
has no solution in integers because of a magnitude issue (i.e. the
left hand side cannot be greater than 3).

This also applies to the equation:
\[
12^{x}+23^{y}=199\hspace{2cm}(x,y\in\mathbb{N}^{0})
\]
because the left hand side is either less than 199 or greater than
199 (i.e. there is no combination of $x,y$ that produces 199 where
we can reach this conclusion easily by just checking the low-values
of $x,y$).

\subsection{Simple Divisibility Checks}

For example, it should be fairly obvious that the equation:
\[
x^{4}-x^{2}+y^{4}-y^{2}=18\hspace{2cm}(x,y\in\mathbb{Z})
\]
has no solution because the left hand side is divisible by 4 $\big[$noting
that $x^{4}-x^{2}=(x^{2}-x)(x^{2}+x)$ where both factors are even
and this similarly applies to $y^{4}-y^{2}$$\big]$ while the right
hand side is not divisible by 4.

\subsection{\label{subSimpleModularCheck}Simple Modular Arithmetic Checks}

For example, it should be fairly obvious (to someone with modest experience
in solving Diophantine equations) that the Diophantine equation:
\[
15x^{2}+6y^{2}=12\hspace{2cm}(x,y\in\mathbb{Z})
\]
has no solution because by a simple modularity inspection (i.e. via
reducing the equation in modulo 5) we find that this equation implies
$y^{2}\stackrel{5}{=}2$ which obviously has no solution (because
2 is a quadratic non-residue of 5) and hence the original Diophantine
equation is not solvable. Also see $\S$ \ref{subModRedPoly} and
$\S$ \ref{secReductionModular}.

It is worth noting that we can consider parity checks (which we investigated
earlier in $\S$ \ref{secParityChecks}) as an instance of simple
modularity (or modular arithmetic) checks. In fact, we can consider
parity checks as the simplest modularity checks (since parity checks
are based on the modular arithmetic of 2 which is the least modulo
in modular arithmetic). However, parity checks (unlike common modularity
checks) are not limited to explicit modularity inspection and checks,
and hence from this perspective we may consider parity checks as more
general than other modularity checks.

\subsection{\label{subFinThoughtInitSensCheck}Final Thought about Initial Sensibility
Checks}

It is important to note that many Diophantine equation problems do
not need in their solution more than an informed inspection based
on these initial and simple checks and hence it is worthwhile to spend
a few minutes on doing this sort of initial inspection and tests which
can save a considerable amount of time and effort in trying to solve
the given problem by the use of sophisticated approaches and techniques
(which may or may not lead to the required result). So in brief, an
initial and systematic inspection using general rules (such as the
rules of parity, primality, sign, divisibility, etc.) can save a lot
of time trying to solve an equation that has no solution or has an
obvious solution and hence it does not require any effort to solve.

It should be obvious that conducting initial sensibility checks as
a first step in dealing with Diophantine problems applies not only
to single Diophantine equations but also to systems of Diophantine
equations. So, if a system contains a non-solvable equation then the
system is not solvable. Therefore, it is worthwhile to inspect the
individual equations of the system (to check if they are solvable
or not) before trying to solve the system. For example, if we are
given the following system:
\[
x^{2}+y^{3}-z^{4}=3\hspace{1.5cm}\textrm{and}\hspace{1.5cm}2x+y^{2}-3y-6z^{2}=75\hspace{2cm}(x,y,z\in\mathbb{Z})
\]
then a quick initial inspection should reveal that the second equation
has no solution (due to parity violation) and hence the system has
no solution (with no need for inspecting the first equation or the
system as a whole).

It is also worthwhile to inspect the characteristics of the system
as a whole (i.e. not only its individual equations separately) to
see if it is sensible for the system to have a solution or not. For
example, if we are given the following system of Diophantine equations:
\[
x^{2}+y^{8}+z^{6}=0\hspace{2cm}\textrm{and}\hspace{2cm}x^{4}+y^{2}-17=0\hspace{2cm}(x,y,z\in\mathbb{Z})
\]
then a quick initial inspection should reveal that the first equation
has only the trivial solution (i.e. $x=y=z=0$), while the second
equation can have only non-trivial solutions (i.e. it cannot accept
the trivial solution). This means that the two equations cannot have
a common solution and hence the system has no solution.

Initial inspection to the system as a whole should also reveal that
the following system (where $x,y,z\in\mathbb{Z}$):
\[
x^{6}+y^{2}-z^{3}=0\hspace{1cm}\textrm{and}\hspace{1cm}x^{2}+y^{4}+z+1=0\hspace{1cm}\textrm{and}\hspace{1cm}x+y+z+1=0
\]
has no solution because if the first equation has a solution then
$z$ must be non-negative, while if the second equation has a solution
then $z$ must be negative. So, the two equations cannot have a common
solution and hence this system cannot have a solution.

\clearpage{}

\section{\label{secPrelimGraphInspect}Graphic Inspection}

Many systems of Diophantine equations can be solved graphically (or
at least graphic investigation of these systems can help in finding
their solutions). This is particularly true when we deal with 2-variable
systems of Diophantine equations of various types (such as polynomials
and exponentials). So, it is recommended to use graphic tools (when
possible and applicable) for initial inspection of systems of Diophantine
equations by plotting the functions representing these equations on
the same graph to see if and where they have points of intersection.
As indicated, the use of these tools can lead to the final solution
of the problem without further action.

In simple cases, ``graphic inspection'' may not require actual plotting
of graphs, i.e. we just ``plot'' the graphs mentally or use graphic
arguments and considerations to reach the final solution. For example,
the following system of Diophantine equations:

\begin{tabular*}{15.95cm}{@{\extracolsep{\fill}}lllll}
\noalign{\vskip0.1cm}
 & $x^{2}-2x+4+y=0$ & and & $x^{2}+y^{2}+6x-10y+30=0$ & \tabularnewline[0.1cm]
\end{tabular*}

\noindent can be solved ``graphically'' by noting that the first
equation represents a parabola while the second equation represents
a circle, and hence the problem can be solved ``graphically'' with
no need to make any plot. So, if we put these equations in their standard
forms (so that we can easily identify the shape and position of the
graphs they represent) then we have:\\
\begin{tabular*}{15.95cm}{@{\extracolsep{\fill}}lllll}
\noalign{\vskip0.1cm}
 & $y+3=-(x-1)^{2}$ & and & $(x+3)^{2}+(y-5)^{2}=4$ & \tabularnewline[0.1cm]
\end{tabular*}\\
Now it is obvious that the first equation represents a parabola which
concaves down with vertex at $(x,y)=(1,-3)$ while the second equation
represents a circle with radius 2 and center at $(x,y)=(-3,5)$ and
hence they cannot have a common point. Therefore, this system obviously
has no solution.

Graphic inspection can also be useful in solving individual Diophantine
equations (i.e. not only systems of such equations) since graphic
inspection can provide an insight about the nature of the equation
and if it can/cannot have a solution (or at least if it has an obvious
solution, e.g. when the graphic plot can show such a solution).

\clearpage{}

\section{\label{secManipulationsTransformation}Manipulations and Transformations}

We should always consider manipulating the given Diophantine equation
(algebraically and non-algebraically) or/and transforming its variables
to put the equation in a ``more friendly'' (or more recognizable
or more familiar) form, e.g. by being of a familiar standard form
whose solution can be obtained more easily or by simplifying its form
and hence making its analysis and investigation more easy. Simple
(algebraic and non-algebraic) manipulations include for instance:

\noindent $\bullet$ Changing the symbols to improve the look of the
equation and remove potential sources of confusion and ambiguity.

\noindent $\bullet$ Moving terms from one side to the other (e.g.
to separate the variables which could make the equation easier to
analyze and tackle).

\noindent $\bullet$ Multiplying/dividing the two sides of the equation
by a constant (e.g. by $-1$ to change signs or by a constant in a
denominator to remove a fraction).

\noindent $\bullet$ Multiplying/dividing the two sides of the equation
by a variable or variables (but we should remember that such manipulations
could introduce or eliminate some solutions and hence all the obtained
solutions after this type of manipulation should be tested on the
original equation as well as considering other potential solutions).

\noindent $\bullet$ Raising to integer powers (e.g. by squaring and
cubing) to remove roots. However, this sort of manipulation can affect
the solutions (e.g. by introducing foreign solutions), and hence we
recommend again testing the obtained solutions after this type of
manipulation on the original equation.

\noindent $\bullet$ Raising to fractional powers which is equivalent
to taking roots (whether integer or fractional roots). Again, this
sort of manipulation can affect the solutions and hence the obtained
solutions after this type of manipulation should be tested on the
original equation.

\noindent $\bullet$ Grouping terms of certain types within brackets
and parentheses. For example, by grouping the terms according to their
variables (e.g. those involving $x$ versus those involving $y$)
the equation can become easier to recognize, categorize and analyze.
Similarly, grouping the terms according to their degree (e.g. linear,
squared, cubic, etc.) can have a similar beneficial effect.

\noindent $\bullet$ Completing squares and hence reducing the size
of the equation and possibly putting it in a certain standard form
whose solution can be obtained more easily. For example, the Diophantine
equation:
\[
x^{2}+y^{2}-6x+12y-36=0
\]
can be put in the following quadratic form:
\[
(x-3)^{2}+(y+6)^{2}=9^{2}
\]
by completing the squares. This new form is much simpler to analyze
than the original equation form and hence it can be solved more easily.
For instance, the new form will naturally indicate the restrictions
$0\leq(x-3)^{2}\leq9^{2}$ and $0\leq(y+6)^{2}\leq9^{2}$ and hence
to find the solutions we need no more than testing the few possibilities
of the values of $x$ and $y$ obtained by considering these restrictions.

We should also consider transforming the variables of the equation
(with or without manipulations which we dealt with already). Transformation
of variables can greatly simplify the given equation and hence make
it easier to analyze and solve. It can also put it in a more recognizable
(and possibly standard) form (such as Pell's equation form; see $\S$
\ref{subPellEquation}) whose solution can be obtained more easily.

For example, with some simple manipulations the equation:
\[
4x^{2}-6y^{2}+12x+108y-478=0
\]

\noindent can be transformed from its (rather messy) form to the tidy
form:
\[
X^{2}-6Y^{2}=1
\]

\noindent where $X=2x+3$ and $Y=y-9$ and hence it can be solved
more easily as a Pell equation where the final solutions are obtained
by the reverse transformations $\big[$i.e. $x=(X-3)/2$ and $y=(Y+9)$$\big]$
subject to certain conditions.

This similarly applies to the equation:
\[
9x^{2}-325y^{2}-42x-130y+35=0
\]

\noindent which can be transformed to the tidy form:
\[
X^{2}-13Y^{2}=1
\]

\noindent where $X=3x-7$ and $Y=5y+1$ (noting that the Pell solutions
will not lead to solutions to the original equation because the reverse
transformations do not produce integer solutions).

We may also employ a single-variable transformation. For example,
the equation:
\[
10x^{2}+2x-8y^{2}=0
\]

\noindent can be transformed to:
\[
X^{2}-80y^{2}=1
\]

\noindent where $X=10x+1$ and hence it can be solved more easily
where the solutions of the original equation are obtained by the reverse
transformation, i.e. $x=(X-1)/10$ subject to certain conditions.

More simple transformations (such as by changing the sign of some
or all variables) may also be used to improve the shape and form of
the equation or to put it in a more solvable form, e.g. by being similar
to an already-solved problem or by being in a certain standard form
whose solution can be obtained readily (see for instance the examples
in $\S$ \ref{secComparisonToSimilar}).

So in brief, manipulations and transformations can generally improve
the ``look and feel'' of the equation (among other beneficial effects)
and hence they usually make it tidier, simpler and easier to analyze
and solve.

\clearpage{}

\section{\label{secInitialComputInv}Initial Computational Investigation}

Initial investigation of Diophantine problems by the use of computational
tools (such as computer codes or spreadsheets or software packages)
should provide an initial impression and insight about the nature
of the expected (and sought-after) solutions. In fact, initial computational
investigation can be a great help in identifying and producing a theoretical
and general argument (or proof or formulation or $\ldots$) that solves
the problem completely and unequivocally.

For example, if we conduct an initial computational investigation
on the Diophantine equation:
\[
4x^{2}+4x-15-y^{3}=0
\]
by a computer code (which loops, for instance, over all integers $x,y$
between $-10000$ and $+10000$) then we should get no solution in
this range and this should suggest that this equation has no solution.
So, if we now inspect this equation further (in the light of this
suggestion) then we may note that completing the square in $x$ will
lead to the equation $y^{3}=(2x+1)^{2}-16$ which can be factorized
as $y^{3}=(2x-3)(2x+5)$. Further reasoning (based on this type of
factorization) should lead to the conclusion that this equation cannot
have a solution.

A similar initial computational investigation on the Diophantine equation:
\[
x^{3}+y^{3}+z^{3}=58
\]
should also fail to find a solution in the given range and this should
suggest that this equation has no solution. Further inspection and
analysis (based, for instance, on modular arithmetic) should produce
a conclusive argument that this equation has no solution (e.g. if
we reduce this equation in modulo 9 then we get the congruence equation:
\[
x^{3}+y^{3}+z^{3}\stackrel{9}{=}4
\]
which has no solution because the sum of three cubes cannot be congruent
to 4 modulo 9 and hence the original Diophantine equation cannot have
a solution).

So, the insight obtained by these initial computational investigations
provides a great help in solving these Diophantine equations (i.e.
by producing a conclusive argument, based on the aforementioned factorization
or modular arithmetic analysis, that these equations have no solution).

Similarly, if we conduct an initial computational investigation on
the following Diophantine equation in a certain range (say as before):
\[
xy+yz=xyz
\]
then we will find out that all the solutions that we get in the given
range have one of 3 forms: $y=0$ and $x$ and $z$ are arbitrary,
$x=z=0$ and $y$ is arbitrary, and $x=z=2$ and $y$ is arbitrary.
Further inspection and analysis (based on the insight provided by
the patterns of the obtained solutions) should lead to a conclusive
argument that this equation actually has only 3 types of solution,
i.e. $(x,y,z)=(m,0,k)$, $(0,n,0)$ and $(2,n,2)$ where $m,n,k\in\mathbb{Z}$.
So, thanks to the insight obtained by this initial and easy computational
investigation we were able to produce a conclusive logical/mathematical
argument that determines all the solutions of this equation unequivocally.
In fact, without such computational investigations solving equations
like this can be much more difficult.

There are many other examples that demonstrate the aid provided by
simple initial computational investigation in tackling Diophantine
equation problems and finding their solutions. Many of the examples
that will be given in the future (for various purposes) are based
on (and initiated by) such computational investigation effort.

As indicated earlier (see the paragraph before the last of $\S$ \ref{secIntroduction}),
computational investigation on its own cannot provide a solution to
Diophantine problems even if this investigation leads to finding all
the actual solutions. So, it should be obvious that ``Initial Computational
Investigation'' (which is the title of this section) is no more than
an initial step in the search for solution (i.e. proving that there
is no solution or finding all the solutions with a conclusive logical/mathematical
argument that there are no other solutions).

\clearpage{}

\section{Classification and Recalling Standard Methods}

The obvious next recommendation (assuming the problem has not been
solved so far) is to classify the problem such as being linear or
non-linear, exponential or polynomial (involving quadratic or/and
cubic or $\ldots$), 2- or 3- $\cdots$ or $n$-variable problem,
and so on. This should obviously be associated with (or followed by)
recalling the standard methods of solution for the specific type (as
identified according to its classification). In this regard, we recommend
using previously-solved similar problems (related to the identified
specific type) as prototypes and models to see if it is possible to
apply their methods of solution to the problem at hand (also see $\S$
\ref{secComparisonToSimilar}). In fact, in some cases solving the
problem at hand may require no more than copying and pasting the solution
of a previously-solved problem with some modifications and adaptations
to reflect the specific characteristics of the problem at hand.

In the following subsections we investigate briefly some common types
of Diophantine problems and some standard methods of their solution.

\subsection{\label{subLinearEquations}Linear Equations}

Solving linear Diophantine equations should be straightforward in
general because there are closed form formulae for their solutions
(at least for the 2- and 3-variable equations). Yes, sometimes the
statement of the problem requires an investigation (or an insight)
about how to model it by a linear Diophantine equation or system of
equations, but this is another story which is related to modeling
rather than solving a modeled (i.e. already-formulated) problem which
is what our paper is about.

In this regard, we should remember that a linear Diophantine equation
has either no solution or infinitely many solutions. However, mathematical
restrictions on the domain or range of solution (e.g. being in natural
numbers, possibly within a given range, instead of being in integers)
should generally reduce the solutions (although it may or may not
affect their infinitude). Physical restrictions on the nature of specific
problem may also reduce the solutions by making the number of solutions
finite or even zero (i.e. when the physical restrictions cannot be
satisfied by the mathematical formulation and hence the equation or
system cannot have a solution within the given physical restrictions).\footnote{The common example of this is when the physical restrictions require
positive values while the mathematical formulation leads to negative
values.}

It is also useful to remember that a linear Diophantine equation is
always solvable (and hence it has infinitely many solutions from a
mathematical viewpoint) when it is homogeneous and when the greatest
common divisor of its coefficients is 1.

Regarding systems of linear Diophantine equations, they generally
can be solved by the well-known methods of linear algebra (with the
restriction on the solutions to be integers) which the interested
reader should refer to in the wide mathematical literature. However,
these systems may also be solved by solving their individual equations
(by the usual techniques and methods of number theory) with taking
the intersection of their solutions (which could be the empty set).\footnote{This issue has been outlined earlier in this paper (see the last paragraph
of $\S$ \ref{secIntroduction}).} For example, if we solve the following system:
\[
15x+10y+30z=41\hspace{1cm}\textrm{and}\hspace{1cm}22x-21y+8z=5\hspace{1cm}\textrm{and}\hspace{1cm}x+19y-39z=73
\]
by the familiar methods of linear algebra we get no integer solution
and hence there is no solution to this system in $\mathbb{Z}$. This
result can also be reached (more simply and directly) by noting that
the left hand side of the first equation is 0 (modulo 5) while its
right hand side is 1 (modulo 5) and hence this equation (as well as
the system) has no solution.

\subsection{\label{subNonLinearPolynomial}Non-Linear Polynomial Equations}

There is no standard method (whether single or multiple) for solving
non-linear polynomial Diophantine equations as they come in many different
shapes and forms (varying, for instance, in the number of variables,
the highest degree of each variable, whether they contain mixed-variables
terms or not, and so on). However, there are a number of recommendations
and considerations (including methods and techniques) that should
be remembered when dealing with non-linear polynomial Diophantine
equations. Some of the most common of these recommendations and considerations
are outlined in the following sub-subsections.

\subsubsection{\label{subFermatLastTheorem}Fermat's Last Theorem}

We should remember Fermat's last theorem when we deal with polynomial
Diophantine equations in 2 or 3 variables with no mixed-variables
terms. This theorem states: no natural numbers $a,b,c$ satisfy the
equation $a^{n}+b^{n}=c^{n}$ for any $n\in\mathbb{N}$ greater than
2. Hence, this theorem is about eliminating the possibility of solutions
of certain type. In this regard it is useful (and possibly important)
to keep the following points in mind when dealing with such equations:
\begin{enumerate}
\item Although this theorem is about ``natural numbers'' it also extends
to negative integers because if $n$ is even then the power does not
distinguish between positive and negative, while if $n$ is odd then
the difference between the negatives and positives is just a multiplicative
factor of $-1$ and hence it does not affect solvability.\footnote{Some cases of mixed positive and negative integers when $n$ is odd
can be dealt with by manipulations (as will be outlined in point \ref{enuP4}).}
\item This theorem is about ``natural numbers'' (as well as their negatives)
and hence the equation $a^{n}+b^{n}=c^{n}$ is generally solvable
when we extend its domain to include 0. In other words, this equation
generally accepts ``trivial'' solutions in the extended sense of
trivial, i.e. when some (and not necessarily all) variables are 0.
This is important to remember to avoid the mistake (which some people
commit) that this equation has only the \textit{trivial} solution
(i.e. $a=b=c=0$) which is not true. For example, $a^{3}+b^{3}=c^{3}$
has solutions like $(a,b,c)=(0,1,1)$ or $(2,0,2)$ which are not
\textit{trivial} although they are ``trivial'' in the aforementioned
extended sense.
\item Fermat's last theorem should also be remembered when dealing with
polynomial Diophantine equations in 2 variables of the above type,
i.e. when dealing with equations of the form $a^{n}+b^{n}=m$ or $a^{n}+b^{n}+m=0$
where $m$ is a given integer (e.g. 27). In other words, we should
inspect $m$ to see if it is an $n^{th}$ power of an integer and
hence the equation can be subject to Fermat's last theorem.
\item \label{enuP4}The above standard form (i.e. $a^{n}+b^{n}=c^{n}$)
of the equation that is subject to Fermat's last theorem can be disguised
and hence we should always inspect the possibility of manipulating
or transforming a ``disguised Fermat's last theorem equation'' to
put it in its standard form (i.e. $a^{n}+b^{n}=c^{n}$). For example,
the following equations are ``disguised Fermat's last theorem equations'':
\[
x^{3}+y^{6}=z^{9}\hspace{1.5cm}x^{3}-y^{3}=z^{3}\hspace{1.5cm}x^{5}-7776=z^{5}\hspace{1.5cm}x^{6}-y^{3}+216z^{3}=0
\]
because these equations can be put in the following standard forms:
\[
x^{3}+Y^{3}=Z^{3}\hspace{1.8cm}y^{3}+z^{3}=x^{3}\hspace{1.8cm}z^{5}+6^{5}=x^{5}\hspace{1.8cm}X^{3}+\zeta^{3}=y^{3}
\]
where $Y=y^{2}$, $Z=z^{3}$, $X=x^{2}$ and $\zeta=6z$. Other forms
of manipulations and transformations (usually more complicated than
the above) can put ``disguised Fermat's last theorem equations''
in their standard form. So, we should always be vigilant and resourceful
(by using manipulations and transformations) so that we can make use
(if possible) of Fermat's last theorem when we deal with equations
of such types.
\end{enumerate}

\subsubsection{\label{subPythagoreanTriples}Pythagorean Triple Rules}

We should also remember the Pythagorean triple rules and theorems
and the equations representing them when dealing with quadratic equations
in 2 or 3 variables with no mixed-variables terms. In this regard,
we should consider manipulating the equation (if necessary) to put
it in a standard Pythagorean triple equation form (e.g. by completing
the squares or/and by transforming the variables; see $\S$ \ref{secManipulationsTransformation}).

For example, if we manipulate the Diophantine equation:
\[
4x^{2}+16y^{2}-z^{2}+4x-24y+10=0
\]
then we get:
\[
(2x+1)^{2}+(4y-3)^{2}=z^{2}
\]
It is fairly obvious that this equation has no solution because ``no
perfect square can be the sum of two odd squares'' which is a rule
based on the Pythagorean triple rules and properties (also see $\S$
\ref{secBasicRules}).\footnote{One of the Pythagorean triple rules is: if $(a,b,c)$ is a primitive
Pythagorean triple then $a$ and $b$ have opposite parity (noting
that multiplication by neither odd factor nor even factor can change
the parity in a way that makes both $a$ and $b$ odd).}

\subsubsection{\label{subPellEquation}Pell's Equation}

We should also remember Pell's equation when we deal with quadratic
polynomial equations in two variables with no mixed-variables terms.
As before, we should consider manipulating or/and transforming the
equation (if necessary) to put it in a standard Pell equation form
(if possible), e.g. by scaling the equation and completing squares.
In fact, such manipulations and transformations are generally useful
(regardless of the possibility of putting the equation in Pell equation
form or not) because they usually simplify the expressions in the
equation and reduce the number of its terms (as well as organizing
the equation in general). All these benefits should help to identify
and recognize potential patterns that can indicate the solutions or
the method of solution (whether by Pell techniques or by something
else).

For example, the following Diophantine equation:
\[
x^{2}-2y^{2}=1\hspace{2cm}(\textrm{or }x^{2}-1=2y^{2})
\]
is already in a standard Pell equation form (or almost) and hence
its solution can be obtained readily by the standard techniques of
solving Pell's equations. So, virtually no effort is required to solve
this equation.

On the other hand, the following Diophantine equation:
\[
4x^{2}-6y^{2}+12x+108y-478=0
\]
 is not in a standard Pell equation form but it can be manipulated
and transformed to a standard Pell equation form with some minimal
effort, that is:
\[
X^{2}-6Y^{2}=1
\]
where $X=2x+3$ and $Y=y-9$. Now, in this form it can be easily solved
as a Pell equation and hence the solutions of the original equation
can be obtained by the reverse transformations, i.e. $x=(X-3)/2$
and $y=Y+9$.

However, it is important to remember the following notes and recommendations
when dealing with polynomial Diophantine equations which are supposed
to be solved (directly or indirectly) by the Pell equation methods:
\begin{enumerate}
\item We should consider in this regard the generalized form of Pell's equation
(i.e. not only its basic form).\footnote{\noindent Pell's equation is generalized to the form $x^{2}-dy^{2}=c$
where $0\neq1\neq c\in\mathbb{Z}$.}
\item When tackling a Pell equation (whether transformed or not) it is recommended
to search for the fundamental solution by inspection before trying
sophisticated processes and techniques (like continued fractions technique).
For example, a few-seconds inspection to the Pell equation:
\[
x^{2}-5y^{2}=1
\]
should lead to the fundamental solution $(x_{1},y_{1})=(9,4)$ and
hence to all other solutions using the well-known formulations and
instructions (related to Pell equation solution). Also basic computational
effort (e.g. through the use of a spreadsheet or a simple code) can
be more economic in the search for the fundamental solution.
\item As soon as we solve a Pell equation obtained by manipulations and
transformations we can try obtaining the solutions of the original
equation by the reverse transformations. However, since we accept
only the integer solutions to the original equation, the existence
of solution to the transformed Pell equation does not guarantee the
existence of solution to the original equation because the reverse
transformations may not produce integer solutions to the original
equation. For example, the following Diophantine equations:
\begin{eqnarray*}
4x^{2}-6y^{2}+12x+108y-478 & = & 0\\
9x^{2}-325y^{2}-42x-130y+35 & = & 0\\
10x^{2}+2x-8y^{2} & = & 0
\end{eqnarray*}
can be easily transformed to the following standard Pell equation
form:
\begin{eqnarray*}
X^{2}-6Y^{2} & = & 1\hspace{2cm}(X=2x+3,\ Y=y-9)\\
X^{2}-13Y^{2} & = & 1\hspace{2cm}(X=3x-7,\ Y=5y+1)\\
X^{2}-80y^{2} & = & 1\hspace{2cm}(X=10x+1)
\end{eqnarray*}
All these Pell equations have solutions. However:\\
$\bullet$ The first original equation accepts all the solutions obtained
by the reverse transformations (because the reverse transformations
produce integer solutions to the original equation in all cases).\\
$\bullet$ The second original equation does not accept any of the
solutions obtained by the reverse transformations (because the reverse
transformations do not produce integer solutions to the original equation
in any case).\\
$\bullet$ The third original equation accepts only some of the solutions
obtained by the reverse transformation (because the reverse transformation
produces integer solutions to the original equation only in some cases).
\item \label{enuPoint4}The usefulness of Pell equation (and the possibility
of its exploitation) in solving Diophantine equations is not limited
to solving the given equation directly (with or without manipulations
and transformations) but it extends beyond this. For example, the
following Diophantine equation:
\[
x^{2}+x-2y^{2}=0
\]
can be treated as a (one-variable) quadratic in $x$ and hence it
has a solution if its discriminant $\Delta$ is a perfect square,
i.e.
\[
\Delta=1^{2}-4(-2y^{2})=k^{2}\hspace{1cm}\rightarrow\hspace{1cm}k^{2}-8y^{2}=1\hspace{2cm}(k\in\mathbb{N})
\]
As the latter is a Pell equation it has solutions (and actually infinitely
many solutions), and hence the original Diophantine equation must
also have solutions (and actually infinitely many solutions).
\end{enumerate}

\subsubsection{\label{subModRedPoly}Modular Reduction}

A well-known approach for solving Diophantine equations in general
(especially those of polynomial types) is the reduction of equation
in an appropriate modulo (and possibly in more than one modulo) which
can reduce the number of variables (or/and show certain features or
patterns) and hence simplify the search for solution. The technique
of modular reduction (which is widely used for solving Diophantine
equations whether of polynomial types or other types) will be investigated
further later on (see $\S$ \ref{secReductionModular}). So, here
we only give some simple examples for the use of this technique for
solving polynomial Diophantine equations.

For example, we can easily \textit{solve} the following polynomial
Diophantine equation (i.e. find out that it has no solution in $\mathbb{Z}$):
\begin{equation}
15x^{2}-35y^{3}=10\label{eqPol15x2}
\end{equation}
by reducing it in modulo 7 to get: $x^{2}\stackrel{7}{=}3$ which
obviously has no solution (since 3 is not a quadratic residue of 7)
and hence we can easily conclude that the original Diophantine equation
has no solution (see $\S$ \ref{subLogFoundModRed}).

Similarly, the Diophantine equation:
\begin{equation}
x^{2}-3y-2z=0\label{eqPolx2}
\end{equation}
can be reduced in modulo 2 to get: $x^{2}-y\stackrel{2}{=}0$ which
is much easier to analyze and solve. Further investigation (based
on considering the parity of $x$ and $y$) should lead to the solution
of $x^{2}-y\stackrel{2}{=}0$ and hence to the solution of $x^{2}-3y-2z=0$.\footnote{The solutions are (where $k,s\in\mathbb{Z}$): $(x,y,z)=(2k,2s,2k^{2}-3s)$
and $(2k+1,2s+1,2k^{2}+2k-3s-1)$.}

More examples about solving non-linear polynomial Diophantine equations
(as well as other types of Diophantine equations) by the use of modular
arithmetic reduction technique will be given in the future.

In this regard, we should remember the following recommendations and
guidelines (or ``rules'') which we can gather from personal experience
(as well as from the literature) about the use of modular reduction
in solving polynomial Diophantine equations:
\begin{enumerate}
\item Always try to use small reduction moduli since they are easier in
management and analysis. In general, when using modular reduction
consider starting your investigation and analysis with the smallest
moduli first (i.e. use moduli in increasing order).
\item In general, modular reduction can lead conclusively to the non-existence
of solution but not to the existence of solution (let alone finding
the specific solution) because having a modular solution is more general
than having a solution to the corresponding Diophantine equation (i.e.
having a modular solution is a necessary but not sufficient condition
for having a solution to the corresponding Diophantine equation).
However, modular reduction can lead to finding the solution (through
initiating or outlining a logical/mathematical argument or by showing
a certain pattern for instance). So, it is useful to try pursuing
the consequences of modular reduction analysis even when modular reduction
lead to having modular solution. Such further analysis can lead to
the conclusion of having a solution to the Diophantine equation (and
even to finding the solution specifically).
\item Consider using more than one reduction modulo (i.e. in more than one
modular reduction operation) associated with comparison and analysis
of the results of the different reduction moduli. Such comparison
and analysis can lead to producing a logical/mathematical argument
that leads to the solution of the problem in hand.
\end{enumerate}

\subsubsection{\label{subFactorizationComparison}Factorization Analysis}

Another well-known approach for solving Diophantine equations in general
(including those of polynomial types) is factorization or/and comparison
of factors. In fact, this method is diverse and versatile and hence
in the following we outline a few common types of this method with
some simple and illuminating examples:
\begin{enumerate}
\item \label{enuFCp1}A common type of this method is to put the given Diophantine
equation in such a form where a factorizable algebraic expression
becomes equal to a (factorizable) specific integer and hence a comparison
between the factors on the two sides of the equation produces systems
of simultaneous equations whose solutions produce all the solutions
of the original Diophantine equation. For example, the following Diophantine
equation:
\[
x^{2}+y^{2}+18xy-x^{2}y^{2}-81=0
\]
can be put in the following form:
\[
(x+y)^{2}-(xy-8)^{2}=17
\]
Now, both sides are factorizable, i.e.
\begin{eqnarray*}
(x+y)^{2}-(xy-8)^{2} & = & (x+y-xy+8)(x+y+xy-8)\\
\textrm{and}\hspace{4cm}17 & = & 1\times17=(-1)\times(-17)
\end{eqnarray*}
Now, if we consider all the possible combinations of these factorization
possibilities of the two sides in both orders then we get four systems
of simultaneous equations whose solutions (if exist) produce all the
solutions of the given Diophantine equation.\\
Other examples of this type are the equations:
\[
5x+xy-2y=0\hspace{2cm}x^{2}-y^{2}+8y-28=0\hspace{2cm}x^{2}-y^{2}-12x-3y+1=0
\]
which can be put in the following forms:
\[
(2-x)(y+5)=10\hspace{1.1cm}(x-y+4)(x+y-4)=12\hspace{1cm}(2x-2y-15)(2x+2y-9)=131
\]
and hence they can be analyzed and solved in a similar manner.
\item \label{enuFCp2}Another type of this method is to manipulate the equation
to make a variable equal to a ratio (or quotient) whose numerator
and denominator can be compared and analyzed to deduce the possible
solutions of the original Diophantine equation. For example, the following
Diophantine equation:
\[
x^{3}y-125x+125=0
\]
can be put in the following form:
\[
y=\frac{125(x-1)}{x^{3}}
\]
Now, it is obvious that $(x-1)$ and $x^{3}$ are coprime and hence
$x^{3}$ must divide 125. This means that $x$ can only be $\pm1$
and $\pm5$ and these values of $x$ should produce all the possible
values of $y$ and hence we get all the possible solutions of the
original Diophantine equation.
\item \label{enuFCp3}An example of a type similar to the previous two types
is the equation:
\[
x^{6}y+xy^{6}-256=0
\]
which can be put in the following two forms:
\[
x(x^{5}y+y^{6})=256\hspace{2cm}\textrm{and}\hspace{2cm}y(x^{6}+xy^{5})=256
\]
These forms should indicate that $x$ and $y$ must be divisors of
256. So, by testing all the possibilities of $x$ being equal to one
of the 18 divisors of 256 and $y$ being equal to one of these 18
divisors we find that only one possibility (i.e. $x=y=2$) satisfies
the given Diophantine equation and hence we obtain the required solution.
\end{enumerate}
Also see $\S$ \ref{secFactorization}.

\subsubsection{\label{subOneVarQuad}One-Variable Quadratic Approach}

In some cases it is possible to solve a quadratic Diophantine equation
in two variables (possibly with a mixed-variable term) by treating
it as a quadratic equation in one variable whose discriminant $\Delta$
can be inspected and analyzed to deduce the solution of the given
Diophantine equation.

For example, the following Diophantine equation:
\[
5x^{2}-8xy+11y^{2}-1175=0
\]
can be treated as a quadratic in $x$ and hence we form and inspect
its discriminant $\Delta$, that is:
\[
\Delta=64y^{2}-20(11y^{2}-1175)=-156y^{2}+23500>0\hspace{1cm}\rightarrow\hspace{1cm}y<\sqrt{\frac{23500}{156}}\simeq12.27
\]
Now, if we consider first only the positive values of $y$ then the
only possibilities we have are $y=1,2,\ldots,12$. On testing these
possible values on the original equation we get only three solutions,
i.e. $(x,y)=(-10,5)$, $(18,5)$ and $(2,11)$. Now, if we note that
the equation is indifferent to change of sign of both variables (noting
the mixed-variable term) then we can conclude that we must have three
other solutions, i.e. $(x,y)=(10,-5)$, $(-18,-5)$ and $(-2,-11)$.
So, by this way we obtain all the possible solutions of this equation
in $\mathbb{Z}$.

A second example is the following Diophantine equation:
\[
x^{2}-xy+6x-y+2=0
\]
which can be put in the following form:
\[
x^{2}+(6-y)x+(2-y)=0
\]
and hence it can be treated as a one-variable quadratic in $x$ whose
discriminant must be a perfect square (say $k^{2}$ where $k\in\mathbb{Z}$),
that is:
\[
\Delta\equiv(6-y)^{2}-4(2-y)=k^{2}\hspace{1cm}\rightarrow\hspace{1cm}(k-y+4)(k+y-4)=12
\]
On inspecting and analyzing this discriminant (using factorization
and comparison as indicated in the last equation; see $\S$ \ref{subFactorizationComparison}),
all the solutions of the given Diophantine equation can be easily
obtained.

A third example is the following Diophantine equation:
\[
x^{2}+xy+2y+2=0
\]
which can be treated as a one-variable quadratic equation in $x$
whose discriminant is:
\[
\Delta\equiv y^{2}-8y-8=k^{2}\hspace{1cm}\rightarrow\hspace{1cm}(y-4-k)(y-4+k)=24
\]
On inspecting and analyzing the last equation (as in the previous
example) we get all the solutions.

Also see the example in point \ref{enuPoint4} of $\S$ \ref{subPellEquation}.

\subsubsection{Separation of Variables with Divisibility Analysis}

Some polynomial Diophantine equations in two variables can be manipulated
such that one variable is separated and equated to a rational fraction
involving functions of only the other variable. For example, the equation:
\[
3x+5xy-6y-5=0
\]
can be manipulated to become:
\[
y=\frac{5-3x}{5x-6}
\]
Now, if we apply a divisibility analysis on the RHS of this equation
then we can conclude that the denominator must be $\pm1$ or $\pm7$.
Further algebraic analysis will lead to the conclusion that the given
Diophantine equation has only one acceptable solution which is $(x,y)=(1,-2)$.

This method similarly applies to the equation:
\[
x^{2}-9x-4y-xy+13=0
\]
which can be manipulated to become:
\[
y=\frac{x^{2}-9x+13}{x+4}=x-13+\frac{65}{x+4}
\]
A divisibility analysis on the RHS will lead to the conclusion that
$x+4$ must be $\pm1$, $\pm5$, $\pm13$, $\pm65$. With further
algebraic analysis all the eight solutions of the given Diophantine
equation can be easily found.

\subsubsection{Other Recommendations and Considerations}

Other recommendations and considerations that should be remembered
when dealing with non-linear polynomial Diophantine equations include
(among many other things):
\begin{enumerate}
\item Considering upper and lower bounds on the potential solutions (see
$\S$ \ref{secUpperLower}).
\item Recalling and employing basic rules and principles (see $\S$ \ref{secBasicRules}).
\item Other general recommendations and guidelines that will be investigated
or outlined later on (see for instance $\S$ \ref{secOtherRecommendations}).
\end{enumerate}

\subsection{\label{subExponentialEquations}Exponential Equations}

Diophantine exponential problems (whether individual equations or
systems of equations) come in many different shapes and forms. Except
in very simple cases, solving these problems generally poses a serious
challenge. In fact, they are usually more difficult to solve than
the more common type of Diophantine polynomial problems. Even searching
for their solutions computationally can pose a challenge (in comparison
for instance to such search for the solutions of polynomial problems)
due to the rapid and explosive rise in the magnitude of the inspected
values which imposes serious limitations on the range of inspected
values unless we have access to exceptional computational facilities
and techniques.

In the following sub-subsections we will try to discuss and investigate
some issues (such as recommendations and techniques) which are useful
to know and remember when tackling simple types of Diophantine exponential
problems. The last sub-subsection of this subsection is dedicated
to some general recommendations and guidelines about the use of the
technique of modular reduction and analysis in solving Diophantine
exponential equations (and equations involving exponentials) noting
that modular reduction and analysis is the most common approach for
tackling and solving Diophantine exponential equations (and equations
involving exponentials) and hence it enjoys a special importance that
requires special attention.

However, before we go through these details we would like to outline
a general strategy for tackling and solving Diophantine exponential
equations. In brief, the solutions of Diophantine exponential equations
(and equations involving exponentials) are usually at the very low
values of their variables (when solutions exist). Hence, the best
strategy for solving these equations is to inspect these equations
computationally (see $\S$ \ref{secInitialComputInv}) for solutions
within a certain range of their variables that includes low values
of their variables (e.g. between 0 and 15) and hence we have three
main cases and scenarios (which determine our strategy and approach):\\
$\bullet$ No solution found: in this case we look for an argument
(e.g. based on parity violation or modular contradiction) to prove
that the given equation has no solution.\\
$\bullet$ A few solutions found (usually in the neighborhood of zero):
in this case we look for an argument showing that there are no other
solutions (i.e. at high values of their variables).\\
$\bullet$ Considerable number of solutions found: in this case we
look for a pattern within the found solutions on which we can build
an argument for the type of solutions (where such an argument should
establish whether there are other solutions or not).

\subsubsection{Magnitude Consideration}

The very basic (and possibly trivial) Diophantine exponential equations
can be solved by just inspecting the limitation on the magnitude of
the two sides of the equation. For example, the following Diophantine
exponential equations (where $x,y\in\mathbb{N}^{0}$):\\
\begin{tabular*}{15.95cm}{@{\extracolsep{\fill}}lllll}
\noalign{\vskip0.1cm}
 & $2^{x}+3^{y}=1$ & $4^{x}+9^{y}=2$ & $5^{x}+7^{y}=40369232$ & \tabularnewline[0.1cm]
\end{tabular*}\\
can be easily ``solved'' by noting that the first equation has no
solution because the LHS cannot be less than 2, and the second equation
has only the trivial solution $x=y=0$ because otherwise the LHS will
be greater than 2, while the third equation (assuming it is solvable)
must have solution(s) only for very low values of its variables and
hence the solution(s) can be obtained computationally (knowing for
sure, by this logical argument, that there is no other solution within
the high values of its variables).\footnote{Actually, the third equation has only one solution which is $x=6$
and $y=9$.}

\subsubsection{\label{subModRedExp}Modular Reduction}

The method of modular reduction is as common in solving exponential
Diophantine equations as in solving polynomial Diophantine equations
(see $\S$ \ref{subModRedPoly} and $\S$ \ref{secReductionModular}).
As indicated earlier (and later), one of the common purposes of modular
reduction is to simplify the equation (e.g. by removing some variables)
or/and showing and revealing some ``hidden'' patterns and clues
that can help in solving the Diophantine equation (noting that the
employed modular reduction may not introduce any ``visible'' change
on the original Diophantine equation apart from converting it from
an ordinary equation to a congruence equation). The logical and mathematical
foundations for the use of modular reduction in solving Diophantine
equations will be outlined later on (see $\S$ \ref{subLogFoundModRed})
and hence the unfamiliar reader should refer to that subsection.

For example, we can easily \textit{solve} the following exponential
Diophantine equation (i.e. find out that it has no solution in $\mathbb{N}^{0}$):
\begin{equation}
4^{x}-12^{y}=a\hspace{2cm}(a=1\textrm{ or }a=6\textrm{ or }a=7)\label{eqExp4x}
\end{equation}
by reducing it in modulo 13 to get: $4^{x}-12^{y}\stackrel{13}{=}a$
which has no solution for $a=1,6,7$ (since $1,6,7$ are not modular
values of $4^{x}-12^{y}$ for any combination of $x,y\in\mathbb{N}^{0}$)
and hence we can easily conclude that the original Diophantine equation
has no solution (see $\S$ \ref{subLogFoundModRed}).

Similarly, the exponential Diophantine equation:
\[
4^{x}-3^{y}=1\hspace{2cm}(x,y\in\mathbb{N}^{0})
\]
can be solved by noting that for $x>1$ we have:
\[
4^{x}-3^{y}=2^{2x}-3^{y}=(2^{3}\times2^{2x-3})-3^{y}=(8\times2^{2x-3})-3^{y}\stackrel{8}{=}-3^{y}\stackrel{8}{\neq}1
\]
and hence it has no solution. So, if it has any solution then $x$
must be 0 or 1. A simple inspection will then reveal that the only
possible solution of this Diophantine equation is $(x,y)=(1,1)$.
The same method and argument apply to the exponential Diophantine
equation:
\[
4^{x}-3^{y}=3\hspace{2cm}(x,y\in\mathbb{N}^{0})
\]
where we can prove that the only possible solution of this equation
is $(x,y)=(1,0)$.

There are many other examples (some of which will be given later)
for the use of modular reduction as a tool for solving exponential
Diophantine equations.

It is important to note that in many cases modular analysis may require
using more than one modulo for reduction and analysis. For example,
the equation:
\begin{equation}
5^{x}-3^{y}=2\label{eqExp5xm3y2}
\end{equation}
can be solved by reducing it in modulo 9 where the solution $(x,y)=(1,1)$
is found for the case of $3^{y}\stackrel{9}{=}3$ (i.e. $5^{x}-3^{1}=2$).
Further analysis in modulo 28 is then applied to show that there is
no solution in the other case (i.e. $3^{y}\stackrel{9}{=}0$ or $x\stackrel{6}{=}5$).

Another example is the equation:
\begin{equation}
3^{x}+5^{y}=z^{2}\label{eqExp3xp5yz2}
\end{equation}
where its analysis in modulo 3 and modulo 9 leads to obtaining its
two solutions, i.e. $(x,y,z)=(1,0,\pm2)$.

Also see the examples of $\S$ \ref{subExpParityModular}.

\subsubsection{\label{subExpParityModular}Parity and Modular Analysis}

In many cases, Diophantine exponential equations can be solved by
a combination of parity and modular analysis (noting that parity analysis
is a form of modular analysis; see $\S$ \ref{subSimpleModularCheck}).
For example, the equation:
\[
2^{x}-3^{y}=1
\]
can be solved by considering the parity of $x$ where the equation
is reduced in modulo 3 for odd $x$ to get the solution $(x,y)=(1,0)$,
and analyzed in modulo 8 for even $x$ to get the solution $(x,y)=(2,1)$.

A similar method applies to the equation:
\[
3^{x}-2^{y}=1
\]
by considering the parity of $x$ where the equation is reduced in
modulo 8 for odd $x$ to get the solution $(x,y)=(1,1)$. The equation
is then analyzed further for even $x$ where it is factorized as $(3^{k}-1)(3^{k}+1)=2^{y}$
and analyzed to infer the power of 2 to which these factors correspond
and this analysis leads to the only solution in this case which is
$(x,y)=(2,3)$.

\subsubsection{\label{subRecomModAnalExp}Recommendations about Modular Analysis}

We would like to draw the attention in this sub-subsection to the
following recommendations and guidelines which we can gather from
our personal experience (as well as from the literature) about the
use of the technique of modular reduction in solving exponential Diophantine
equations:\footnote{In fact, most of these recommendations and guidelines also apply to
other types of Diophantine equations.}
\begin{enumerate}
\item Try to find and use a reduction modulo that is as small as possible
so that you deal with small or modest residue system that is easy
to manage and deal with. A systematic inspection to candidate reduction
moduli in their increasing magnitude should lead to the identification
of the smallest modulo that is appropriate for that purpose (if it
exists).
\item Try to choose (when you have a choice) a reduction modulo $k$ so
that the base of interest $n$ for that reduction modulo has a low
integer order, i.e. $O_{k}n$.\footnote{We note that ``the base of interest'' is more general than dealing
with only one powered base (i.e. when modular reduction eliminates
the powers of other bases) and dealing with more than one powered
base (i.e. when more than one powered base remain after modular reduction).} For instance, if the base of interest is 7 (e.g. we are dealing with
$7^{x}$) then $O_{10}7=4$ while $O_{13}7=12$ and hence if we have
a choice then 10 is a better reduction modulo. The reason for this
is that low integer order means fewer cases to deal with and easier
investigation and management.
\item Try to choose (when you have a choice) a reduction modulo that is
prime (not composite) or at least a power of prime. This should reduce
the complications due to the obvious advantages of dealing with prime
moduli (although this can be in conflict with the previous recommendation).
\item Consider using more than one reduction modulo (i.e. in more than one
modular reduction operation) associated with comparison and analysis
of the results of the different reduction moduli. Such comparison
and analysis can lead to producing a logical/mathematical argument
that leads to the solution of the problem in hand. Some examples for
using more than one modulo were given earlier (see for instance the
analysis of Eqs. \ref{eqExp5xm3y2} and \ref{eqExp3xp5yz2}).
\item There is no guarantee that the problem in hand can be tackled by modular
reduction in an appropriate way (commensurate with the nature and
size of the problem) and hence be prepared for possible failure. Accordingly,
be prepared to give up after reasonable time and effort (i.e. don't
be bogged down!).\footnote{This recommendation (in our view) is especially important when dealing
with exponential equations because the search for an appropriate reduction
modulo in this case is usually demanding and time consuming and hence
it can waste considerable amount of time if we keep trying and trying.}
\end{enumerate}

\subsection{Mixed Polynomial-Exponential Equations}

There is no single standard method or technique for solving this type
of Diophantine equations. However, the methods and techniques used
for solving polynomial and exponential equations (see $\S$ \ref{subLinearEquations},
\ref{subNonLinearPolynomial} and \ref{subExponentialEquations})
are generally used to solve this type of Diophantine equations (since
the equation is a mix of these types of Diophantine equations). It
is important to note in this context the special importance of using
modular reduction which is seemingly the most common technique in
tackling this type of problems (although it is usually associated
with other techniques and tricks and further modular and non-modular
analysis). This should be understandable noting the ability of modular
reduction to eliminate one type of the mix (i.e. the polynomial type
and the exponential type) if such elimination is required as well
as its exceptional ability to simplify the equation in general (with
and without elimination of one type) and expose it to the versatile
(and relatively simple) tools and rules of modular arithmetic. We
should also mention in this regard that the use of modular reduction
is common in tackling both polynomial and exponential types of Diophantine
equations (see $\S$ \ref{subModRedPoly} and \ref{subModRedExp})
and hence it should also be common in tackling mixed polynomial-exponential
Diophantine equations.

In the following sub-subsections we will present a sample of this
type of Diophantine equations and the methods used in their solutions.
This should give an idea about the methods used in solving this type
of equations.

\subsubsection{Simple Inspection}

Some mixed polynomial-exponential Diophantine equations are so trivial
that they can be ``solved'' by just simple inspection. For example,
if we write the equation:
\[
7^{x}-8^{y}-z=0\hspace{2cm}(x,y\in\mathbb{N}^{0},\ z\in\mathbb{Z})
\]
as $z=7^{x}-8^{y}$ then it is obvious that there is no restriction
on $z$ (other than being integer) and hence we are free to assign
any value (within the domain) to its two exponential variables which
they fix (as soon as they are assigned) the value of $z$. So, the
solutions of this equation are obviously $(x,y,z)=(k,s,7^{k}-8^{s})$
where $k,s\in\mathbb{N}^{0}$.

\subsubsection{\label{subInspectModularAna}Inspection with Modular Analysis}

For example, the equation:
\begin{equation}
2^{x}+3^{y}=z^{2}\hspace{2cm}(x,y\in\mathbb{N}^{0},z\in\mathbb{Z})\label{eq2xp2yEqz2}
\end{equation}
can be inspected (aided by modular analysis) for low values of $x$
where six solutions can be found, i.e. $(x,y,z)=(0,1,\pm2)$, $(3,0,\pm3)$
and $(4,2,\pm5)$. On applying further modular analysis (in modulo
8) we can conclude that this Diophantine equation has no solution
other than these six. This method also applies to the equation:
\begin{equation}
2^{x}-3^{y}=z^{2}\hspace{2cm}(x,y\in\mathbb{N}^{0},z\in\mathbb{Z})\label{eq2xm3y}
\end{equation}
whose five solutions at low values of $x$ $\big[$namely $(x,y,z)=(0,0,0)$,
$(1,0,\pm1)$ and $(2,1,\pm1)$$\big]$ can be easily found by inspection
while modular analysis (in modulo 8) can be used to prove that it
has no other solutions.

Another example of this method is the equation:
\[
3^{x}-4^{y}=z^{2}\hspace{2cm}(x,y\in\mathbb{N}^{0},z\in\mathbb{Z})
\]
which can be inspected for $x=0$ to find the trivial solution $(x,y,z)=(0,0,0)$.
Modular analysis (in modulo 3) can then easily reveal that this equation
has no solution for $x>0$ and hence the equation has only the trivial
solution. This method also applies to the equation:
\[
4^{y}-3^{x}=z^{2}\hspace{2cm}(x,y\in\mathbb{N}^{0},z\in\mathbb{Z})
\]
whose three solutions at low values of $y$ $\big[$namely $(x,y,z)=(0,0,0)$,
($1,1,\pm1)$$\big]$ can be easily found by inspection while modular
analysis (in modulo 8) will reveal that it has no other solutions.

\subsubsection{\label{subComparison}Comparison to a Similar Equation}

For example, the equation:
\[
4^{x}+3^{y}=z^{2}\hspace{2cm}(x,y\in\mathbb{N}^{0},z\in\mathbb{Z})
\]
can be solved by writing it as $2^{X}+3^{y}=z^{2}$ (where $X=2x$)
and hence comparing it to Eq. \ref{eq2xp2yEqz2}. As we see, this
equation is the same as Eq. \ref{eq2xp2yEqz2} (with $X$ replacing
$x$). Now, if we note that $X=2x$ then we can conclude (by using
the solutions of Eq. \ref{eq2xp2yEqz2}) that the only solutions to
the given equation are $(x,y,z)=(0,1,\pm2)$ and $(2,2,\pm5)$.

\subsubsection{Modular Reduction with Further Analysis}

For example, the equation:
\[
5^{x}-6y+21=0\hspace{2cm}(x\in\mathbb{N}^{0},y\in\mathbb{Z})
\]
can be solved by modular reduction (mod 6) to get $5^{x}+21\stackrel{6}{=}0$,
i.e. $(-1)^{x}+3\stackrel{6}{=}0$. This congruence equation has no
solution (because the left hand side is either 2 or 4) and hence the
given Diophantine equation has no solution.

Similarly, the equation:
\[
3^{x}+5^{y}-4z-2=0\hspace{2cm}(x,y\in\mathbb{N}^{0},z\in\mathbb{Z})
\]
can be solved by modular reduction (mod 4) to get $3^{x}+5^{y}-2\stackrel{4}{=}0$,
i.e. $(-1)^{x}+(1)^{y}-2\stackrel{4}{=}0$. The solution of this equation
is obviously $x=2k$ and $y=s$ ($k,s\in\mathbb{N}^{0}$). On solving
the given equation for $z$ we get $z=(3^{2k}+5^{s}-2)/4$ and hence
the solutions of the given equation are all triples of the following
form: $(x,y,z)=\left(2k,\,s,\,\frac{3^{2k}+5^{s}-2}{4}\right)$. It
is worth noting that $(3^{2k}+5^{s}-2)/4$ is integer for all $k,s\in\mathbb{N}^{0}$
(since $3^{2k}+5^{s}-2$ is zero in modulo 4) and hence this solution
applies to all $k,s\in\mathbb{N}^{0}$.

\subsubsection{Modular Reduction with Substitution}

For example, the equation:
\[
5x+4^{y}-11=0\hspace{2cm}(x\in\mathbb{Z},y\in\mathbb{N}^{0})
\]
can be solved by modular reduction (mod 5) to get $4^{y}-11\stackrel{5}{=}0$,
i.e. $(-1)^{y}-1\stackrel{5}{=}0$. The solution of this congruence
equation is all even $y\geq0$, i.e. $y=2k$ ($k\in\mathbb{N}^{0}$).
On substituting this into the given Diophantine equation and solving
the resulting equation for $x$ we get $x=(11-4^{2k})/5$ and hence
the solutions of the given equation are all pairs of the following
form: $(x,y)=\left(\frac{11-4^{2k}}{5},\,2k\right)$ where $k\in\mathbb{N}^{0}$.
It is worth noting that $(11-4^{2k})/5$ is integer for all values
of $k$ because for $k=0$ it is equal to 2, while for $k>0$ the
numerator $(11-4^{2k})$ ends in 5 and hence it is divisible by 5.

\subsubsection{Modular Reduction with Parity Analysis}

For example, the equation:
\[
5^{x}-11x+3y+1=0\hspace{2cm}(x\in\mathbb{N}^{0},y\in\mathbb{Z})
\]
can be solved by modular reduction (mod 3) to get $(-1)^{x}+x+1\stackrel{3}{=}0$
whose solutions (which can be inferred by parity analysis of $x$
aided by modular inspection) are $x=3+6k$ and $x=4+6k$ where $k\in\mathbb{N}^{0}$.\footnote{We note that if we reduce the equation (mod 3) to the form $5^{x}+x+1\stackrel{3}{=}0$
then we can use a different method of analysis (based on inspection
and induction).} On substituting these expressions of $x$ in the original Diophantine
equation and solving for $y$ we get the required solutions, i.e.
\[
(x,y)=\left(3+6k,\,\frac{-5^{3+6k}+32+66k}{3}\right)\hspace{0.5cm}\textrm{and}\hspace{0.5cm}(x,y)=\left(4+6k,\,\frac{-5^{4+6k}+43+66k}{3}\right)
\]
It is straightforward to show that $(-5^{3+6k}+32+66k)/3$ and $(-5^{4+6k}+43+66k)/3$
are always integers and hence these solutions are valid for all values
of $k\in\mathbb{N}^{0}$.

\subsubsection{Classification with Parity Analysis}

For example, the equation:
\[
7^{x}-8^{y}-2z=0\hspace{2cm}(x,y\in\mathbb{N}^{0},\ z\in\mathbb{Z})
\]
can be solved by classifying it according to the value of its exponential
variables associated with parity analysis. More specifically, the
case $x=0$ and $y>0$ as well as the case $x>0$ and $y>0$ have
no solution due to parity violation. So, all we need to consider is
the case $x=y=0$ which leads to the trivial solution $(x,y,z)=(0,0,0)$
and the case $x>0$ and $y=0$ which leads to the obvious solution
$(x,y,z)=\left(k,0,\frac{7^{k}-1}{2}\right)$ where $k\in\mathbb{N}$.

\subsection{Equations Involving Roots}

Equations involving roots may not be classified technically as Diophantine
equations although this will not prevent us from including them in
our investigation due to the obvious merit and justification of this
inclusion and their undeniable qualification to be treated as such.
There are a variety of methods (or techniques or tricks or $\ldots$)
for tackling Diophantine equations involving roots. In the following
sub-subsections we present some of the most common of these methods
with illuminating and illustrating examples.

\subsubsection{Simple Inspection}

Some of the Diophantine equations involving roots are so trivial that
they can be ``solved'' by just simple inspection. For example, the
equation:
\[
\sqrt{x}-y=379
\]
can be solved by noting that $x$ must be a perfect square (i.e. $x=s^{2}$
where $s\in\mathbb{Z}$) and hence $y=\sqrt{x}-379=|s|-379$ (where
$|s|$ is the absolute value of $s$). So, the solutions are $(x,y)=(s^{2},|s|-379)$
where $s\in\mathbb{Z}$.

Another example is the equation:
\[
5^{x}-7^{y}-2\sqrt{z}=0\hspace{2cm}(x,y,z\in\mathbb{N}^{0})
\]
which can be written as $\sqrt{z}=(5^{x}-7^{y})/2$ and hence $z=\left(\left[5^{x}-7^{y}\right]/2\right)^{2}$
where $(5^{x}-7^{y})\geq0$. So, the solutions of the given equation
are all triples of the following form: $(x,y,z)=\left(k,\,s,\,\left[\frac{5^{k}-7^{s}}{2}\right]^{2}\right)$
where $k,s\in\mathbb{N}^{0}$ and $(5^{k}-7^{s})\geq0$. It is worth
noting that $z$ is always integer because $(5^{k}-7^{s})$ is even
for all $k,s\in\mathbb{N}^{0}$.

\subsubsection{Enumeration}

For example, the equation:
\[
\sqrt{x}+\sqrt{y}=9
\]
can be solved by noting that there are only 10 pairs of $\sqrt{x}$
and $\sqrt{y}$ that can add up to 9, i.e. $(\sqrt{x},\sqrt{y})=(0,9),(1,8)\ldots(9,0)$.
So, if we square the numbers in each pair then we get all the possible
solutions. So, the 10 solutions are: $(x,y)=(0,81)$, $(1,64),\ldots(81,0)$.

This similarly applies to the equation:
\[
\sqrt{x}+\sqrt{y}=\sqrt{363}
\]
which can be solved by noting that $\sqrt{363}=11\sqrt{3}$ and hence
we have only 12 pairs of $\sqrt{x}$ and $\sqrt{y}$ that can add
up to $\sqrt{363}$, i.e. $(\sqrt{x},\sqrt{y})=(0\sqrt{3},11\sqrt{3}),(1\sqrt{3},10\sqrt{3}),\ldots,(11\sqrt{3},0\sqrt{3})$.
So, if we square the numbers in each pair then we get all the 12 possible
solutions, i.e. $(x,y)=(0,363)$, $(3,300),\ldots(363,0)$.

\subsubsection{Linearization}

For example, the equation:
\[
6x+10\sqrt{y}-19z=0
\]
can be solved by letting $Y=\sqrt{y}$ and hence we get $6x+10Y-19z=0$
whose solution is $(x,Y,z)=(11s-19k,\ s,\ 4s-6k)$ where $s,k\in\mathbb{Z}$.
However, $\sqrt{y}$ must be an integer and hence we must have $y=t^{2}$
($t\in\mathbb{Z}$). Therefore, $Y=\sqrt{y}=\sqrt{t^{2}}=|t|=s$.
Thus, the solution of the original equation is $(x,y,z)=(11|t|-19k,\ t^{2},\ 4|t|-6k)$
where $t,k\in\mathbb{Z}$.

Another example of linearization is the equation:
\[
21x+35\sqrt{y}-12\sqrt{z}=41
\]
which can be solved by letting $Y=\sqrt{y}$ and $Z=\sqrt{z}$ and
hence we get $21x+35Y-12Z=41$ whose solution is $(x,Y,Z)=(2-5s-4k,\ 1+3s,\ 3-7k)$
where $s,k\in\mathbb{Z}$. Now, $Y=\sqrt{y}\geq0$ and hence $s$
must be $\geq0$ (i.e. $s\in\mathbb{N}^{0}$). similarly, $Z=\sqrt{z}\geq0$
and hence $k$ must be $\leq0$ (i.e. $\mathbb{Z}\ni k\leq0$). Therefore,
the solution of the original equation is $(x,y,z)=(2-5s-4k,\ [1+3s]^{2},\ [3-7k]^{2})$
where $s\in\mathbb{N}^{0}$ and $\mathbb{Z}\ni k\leq0$.

The equation:
\[
\sqrt{x+1}-\sqrt{y+5}=1
\]
can also be solved by linearization by letting $X=(x+1)$ and $Y=(y+5)$
and hence we have $\sqrt{X}-\sqrt{Y}=1$. Now, if $X=k^{2}$ ($k\in\mathbb{Z}$)
then $\sqrt{Y}=\sqrt{X}-1=|k|-1$, i.e. $Y=(|k|-1)^{2}$. However,
since $\sqrt{X}-\sqrt{Y}=1$ then we must have $X>Y$ which implies
$|k|>(1/2)$, i.e. $k\in\mathbb{Z}$ and $k\neq0$. Hence, $X=(x+1)=k^{2}$
and $Y=(y+5)=(|k|-1)^{2}$ where $k\in\mathbb{Z}$ and $k\neq0$.
So, the solutions of the given equation are all pairs of the following
form: $(x,y)=\left(k^{2}-1,\,\{|k|-1\}^{2}-5\right)$ where $k\in\mathbb{Z}$
and $k\neq0$.

\subsection{\label{subEquationsInvolvingFractions}Equations Involving Fractions}

Again, equations involving fractions may not be classified technically
as Diophantine equations although this will not affect our decision
to include them in our investigation (justified at least by practical
considerations even if they are not Diophantine equations from a theoretical
viewpoint due to the limitation of the definition of Diophantine equations).
Anyway, there are various methods and techniques for tackling Diophantine
equations involving fractions (depending for instance on their types
and number of variables). In the following sub-subsections we will
outline some of the most common of these methods presented and demonstrated
within illuminating examples.

\subsubsection{Magnitude Analysis}

For example, the equations:
\[
\frac{1}{x}+\frac{1}{y}+\frac{1}{z}=4\hspace{2cm}\textrm{and}\hspace{2cm}\frac{1}{x}+\frac{1}{y}+\frac{1}{z}=3\hspace{2cm}(x,y,z\in\mathbb{Z},\ xyz\neq0)
\]
can be easily solved by magnitude analysis where the first equation
has no solution because its LHS cannot be greater than 3 (i.e. when
$x=y=z=1$), while the second equation has only one solution (i.e.
$x=y=z=1$) because if any one of the variables is not equal to 1
then the LHS will be less than 3.

\subsubsection{Sign and Magnitude Analysis}

For example, the equation:
\[
\frac{x}{y}+\frac{y}{x}=1\hspace{2cm}(x,y\in\mathbb{Z},\ xy\neq0)
\]
can be solved by noting that if this equation has a solution then
$x$ and $y$ must have the same sign (because otherwise the sum will
be negative) and hence (whether $x\geq y$ or $x<y$) the LHS must
be greater than 1, i.e. the equation has no solution.

Many similar equations can be analyzed and solved by this method,
e.g. the equation:
\[
\frac{x}{y}+xy=1\hspace{2cm}(x,y\in\mathbb{Z},\ y\neq0)
\]
has no solution because $x$ and $y$ must have the same sign (otherwise
the sum will be negative) and hence (whether $x\geq y$ or $x<y$)
the LHS must be greater than 1 (noting that $x=0$ is not a possibility).

This also applies to the equations:
\[
\frac{1}{x}+\frac{1}{y}=\frac{2}{3}\hspace{2cm}\frac{1}{x}+\frac{2}{y}=\frac{3}{4}\hspace{2cm}(x,y\in\mathbb{Z},\ xy\neq0)
\]
which can be solved by this method (where the details can be found
in \cite{SochiBookNTV2}).

\subsubsection{Divisibility Analysis}

Some equations can be manipulated in one (or more) form that facilitates
simple divisibility analysis which leads to the solutions. For example,
the equation:
\[
\frac{1}{x}+\frac{1}{y}=z\hspace{2cm}(x,y,z\in\mathbb{Z},\ xy\ne0)
\]
can be manipulated into the following two forms:
\[
1+\frac{x}{y}=xz\hspace{3cm}\frac{y}{x}+1=yz
\]
which imply $y=\pm x$ and this (with extra basic analysis) will lead
to the required solutions.

\subsubsection{Separation of Variables with Divisibility Analysis}

For example, the equation:
\begin{equation}
\frac{14}{x}+\frac{y}{19}=25\hspace{2cm}(x,y\in\mathbb{Z},\ x\neq0)\label{eqFrac14x}
\end{equation}
can be manipulated to become:
\[
y=475-\frac{266}{x}
\]
Now, by a basic divisibility analysis (i.e. $x$ must be a divisor
of 266) we can easily obtain all the (sixteen) solutions of the original
Diophantine equation.

This method also applies to the equation:
\begin{equation}
\frac{20}{x}+\frac{33}{y}=2\hspace{2cm}(x,y\in\mathbb{Z},\ xy\neq0)\label{eqFrac20x}
\end{equation}
which can be manipulated to become:
\[
x=\frac{20y}{2y-33}
\]
where by a simple divisibility analysis (i.e. $2y-33$ must be a divisor
of 330) we can obtain all the (fifteen) solutions of the given Diophantine
equation.

Another example of this approach is the equation:
\[
\frac{x}{8}+\frac{y}{5}-\frac{3}{z}=7\hspace{2cm}(x,y,z\in\mathbb{Z},\ z\ne0)
\]
which can be reduced to the form:
\[
z=\frac{120}{5x+8y-280}
\]
and hence it is solved by divisibility analysis.

\subsubsection{Conversion to Polynomial Equation}

Some Diophantine equations involving fractions can be converted to
polynomial equations by multiplying the entire equation by a suitable
factor (and hence it is solved as a polynomial Diophantine equation;
see $\S$ \ref{subLinearEquations} and $\S$ \ref{subNonLinearPolynomial}).
For example, the equation:
\begin{equation}
\frac{x}{y}+\frac{y}{x}=2\hspace{2cm}(x,y\in\mathbb{Z},\ xy\neq0)\label{eqxyPyx2}
\end{equation}
can be solved by multiplying the equation by $xy$ and rearranging
to get $x^{2}+y^{2}-2xy=0$, i.e. $(x-y)^{2}=0$. Hence, we conclude
that the general solution of the given equation is $(x,y)=(k,k)$
where $k\in\mathbb{Z}$ and $k\neq0$.

In fact, this method applies to the more general version of Eq. \ref{eqxyPyx2},
i.e.
\[
\frac{x}{y}+\frac{y}{x}=z\hspace{2cm}(x,y,z\in\mathbb{Z},\ xy\neq0)
\]
to prove that this equation has no solution except for $z=\pm2$.
This is achieved by forming the discriminant of the equation $x^{2}-zxy+y^{2}=0$
(treated as a quadratic in $x$ or as a quadratic in $y$; see $\S$
\ref{subOneVarQuad}) where a detailed analysis of the discriminant
will lead to the required conclusion (i.e. there is no solution except
for $z=\pm2$) as well as finding the solutions when $z=\pm2$.

Another example is the equation:
\begin{equation}
\frac{x}{y}-\frac{y}{x}=1\hspace{2cm}(x,y\in\mathbb{Z},\ xy\neq0)\label{eqFracxyMyx1}
\end{equation}
which can be solved by multiplying it by $xy$ and rearranging to
get $x^{2}-xy-y^{2}=0$. On treating this as a one-variable quadratic
equation in $x$ (see $\S$ \ref{subOneVarQuad}) and analyzing its
discriminant we conclude that the given equation has no solution.

A similar method applies to the equation:
\[
\frac{x}{y}-\frac{y}{x}=2\hspace{2cm}(x,y\in\mathbb{Z},\ xy\neq0)
\]
which can be solved by multiplying it by $xy$ and rearranging to
get $x^{2}=2xy+y^{2}$ and hence:
\[
2x^{2}=x^{2}+x^{2}=x^{2}+(2xy+y^{2})=x^{2}+2xy+y^{2}=(x+y)^{2}
\]
On taking the square root of both sides we get: $x\sqrt{2}=\pm(x+y)$
which is impossible because $x\sqrt{2}$ is irrational while $(x+y)$
is an integer. So, we conclude that the given equation has no solution.

In fact, the method of Eq. \ref{eqFracxyMyx1} applies to the more
general version of Eq. \ref{eqFracxyMyx1}, i.e.
\[
\frac{x}{y}-\frac{y}{x}=z\hspace{2cm}(x,y,z\in\mathbb{Z},\ xy\neq0)
\]
to prove that this equation has no solution except for $z=0$ where
this is achieved by analyzing the discriminant of the equation $y^{2}+zxy-x^{2}=0$
(which is obtained by multiplying the given equation by $xy$) as
a quadratic in $y$ (see $\S$ \ref{subOneVarQuad}). The analysis
will lead to the required conclusion as well as finding the solutions
when $z=0$.

Conversion to polynomial equation can also be used to solve the equation:
\[
\frac{1}{x}+\frac{1}{y}=1\hspace{2cm}(x,y\in\mathbb{Z},\ xy\neq0)
\]
where it is reduced to the form $(1-x)(1-y)=1=(-1)(-1)$ and hence
its only solution (i.e. $x=y=2$) can be obtained by factorization
analysis (as indicated).

\subsubsection{\label{subCompareFrac}Comparison to a Similar Equation}

Some equations can be easily solved by comparing them to similar equations
whose solutions are known (or can be obtained easily). For example,
the equation:
\[
\frac{14}{x}-\frac{y}{19}=25\hspace{2cm}(x,y\in\mathbb{Z},\ x\neq0)
\]
can be easily solved by comparing it to Eq. \ref{eqFrac14x} where
we write our equation as $\frac{14}{x}+\frac{Y}{19}=25$ (with $Y=-y$)
and hence we should have the same solutions for ($x,Y$) as those
found for Eq. \ref{eqFrac14x}. The final solutions for our equation
are then obtained by reversing the sign of $Y$.

This method of comparison similarly applies to the equation:
\begin{equation}
\frac{20}{x}-\frac{33}{y}=2\hspace{2cm}(x,y\in\mathbb{Z},\ xy\neq0)\label{eqFrac20x2}
\end{equation}
which can be compared to Eq. \ref{eqFrac20x} and hence its solutions
are obtained from the solutions of Eq. \ref{eqFrac20x} by reversing
the sign of $y$. This is because if we write Eq. \ref{eqFrac20x2}
as $\frac{20}{x}+\frac{33}{-y}=2$ then it is no more than Eq. \ref{eqFrac20x}
with the sign of $y$ being reversed.

\subsubsection{Symmetry-Based Ordering and Magnitude Analysis}

Some equations have a symmetry in their variables that can be exploited
in an analysis based on ordering the magnitude of the variables. For
example, the equations:
\[
\frac{1}{x}+\frac{1}{y}+\frac{1}{z}=1\hspace{2cm}\textrm{and}\hspace{2cm}\frac{1}{x}+\frac{1}{y}+\frac{1}{z}=2\hspace{2cm}(x,y,z\in\mathbb{N})
\]
can be investigated preliminarily under the ordering assumption that
$x\leq y\leq z$ and hence their solutions (under this assumption)
are obtained rather easily by exploiting the restrictions on the magnitude
of the variables where these restrictions are based on the restriction
on the magnitude of the LHS of these equations (which is imposed by
their RHS) as well as the restriction imposed by the ordering assumption.
The remaining solutions are then obtained by lifting the ordering
assumption and hence permuting the variables in the \textit{initial
solutions} (i.e. the solutions obtained under the ordering assumption)
where permuting the variables in the initial solutions is justified
by the aforementioned symmetry. Hence, all the solutions are obtained
in this two-stage process that exploits and employs symmetry, ordering
and magnitude restrictions.

Also see $\S$ \ref{secSymmetriesCycling}.

\subsection{Equations Involving Roots and Fractions}

The best approach for tackling this sort of equations is to linearize
the radicals (i.e. the roots are replaced by non-radical symbols such
as replacing $\sqrt{x}$ by $X$) and hence the problem is reduced
in difficulty because we deal first with an equation that involves
non-radical fractions while dealing with the issue of radicals is
deferred to the end (where it is usually managed rather easily by
raising the variables of the obtained solutions to suitable powers).
For example, the equation:
\begin{equation}
\frac{1}{\sqrt{x}}+\frac{1}{\sqrt{y}}+\frac{1}{\sqrt{z}}=2\hspace{2cm}(x,y,z\in\mathbb{N})\label{eqRootFracxyz}
\end{equation}
can be solved by converting it first to the form:
\begin{equation}
\frac{1}{X}+\frac{1}{Y}+\frac{1}{Z}=2\hspace{2cm}(X,Y,Z\in\mathbb{N})\label{eqRootFracXYZ}
\end{equation}
where the latter equation is solved as an equation involving fractions
(see $\S$ \ref{subEquationsInvolvingFractions}). The final solutions
(i.e. the solutions of Eq. \ref{eqRootFracxyz}) are then obtained
by squaring the values of the variables of the solutions of Eq. \ref{eqRootFracXYZ}.
Many other similar equations can be dealt with by the same method.

\subsection{Equations involving Factorials}

A common approach in this type of problems is to obtain the solutions
for the low values of factorial by inspection while the cases of high
values of factorial are analyzed and tackled by a more general approach
(e.g. by modular reduction). However, other (more versatile) tricks
and techniques may be required or employed to solve this type of problems
(e.g. when the factorial term is scaled by an integer factor). Some
of these issues and details are demonstrated and clarified in the
following examples:
\begin{enumerate}
\item The equation:
\[
11x-y!=13\hspace{2cm}(x\in\mathbb{Z},\ y\in\mathbb{N}^{0})
\]
has no solution for $y<13$ (where this result is obtained by inspection).
For $y\geq13$ we have $11x\stackrel{13}{=}0$ whose solution is $x=13k$
where $k\in\mathbb{Z}$. On substituting this in the given equation
and simplifying we get $k=\frac{(y!/13)+1}{11}$ which implies that
$k$ cannot be an integer and hence the given equation has no solution.
\item The equation:
\[
x^{2}-y!=3\hspace{2cm}(x\in\mathbb{Z},\ y\in\mathbb{N}^{0})
\]
has no solution for $y>3$ because $x^{2}\stackrel{4}{=}3$ (which
we obtain by reducing the given Diophantine equation in modulo 4)
has no solution (since 3 is not a quadratic residue of 4). Hence,
if there is any solution then we must have $y=0,1,2,3$. So, by inspecting
these values of $y$ we get all the solutions of the given equation,
i.e. $(x,y)=(\pm2,0)$, $(\pm2,1)$, $(\pm3,3)$.\\
The equation:
\[
x^{2}-y!=2\hspace{2cm}(x\in\mathbb{Z},\ y\in\mathbb{N}^{0})
\]
can be similarly analyzed and solved (where reduction in modulo 3
is used for $y>2$ with inspection of $y=0,1,2$).
\item The equation:
\[
16x+9y!=121\hspace{2cm}(x\in\mathbb{Z},\ y\in\mathbb{N}^{0})
\]
has no solution for $y>1$ (due to parity violation) and hence it
is solvable only for $y=0,1$ (i.e. $x=7$ corresponding to $y=0,1$).
\item The equation:
\[
x^{2}+y^{2}-z!=3\hspace{2cm}(x,y\in\mathbb{Z},\ z\in\mathbb{N}^{0})
\]
can be solved by noting that for $z>3$ we have $x^{2}+y^{2}\stackrel{4}{=}3$
which has no solution and hence all we need to do is to inspect $z=0,1,2,3$
which lead to its (twenty) solutions.
\item The equation:
\[
x^{2}+y^{2}+z!=24\hspace{2cm}(x,y\in\mathbb{Z},\ z\in\mathbb{N}^{0})
\]
can be solved by noting that $z$ cannot be greater than 4 and hence
all we need to do is to inspect $z=0,1,2,3,4$ which lead to its (five)
solutions.
\item The equation:
\[
3x+5y+15z!=17\hspace{2cm}(x,y\in\mathbb{Z},\ z\in\mathbb{N}^{0})
\]
can be solved by reducing it in modulo 5 to get $x=4+5k$ and by reducing
it in modulo 3 to get $y=1+3s$. On substituting these into the original
equation we get its (parameterized) solutions.
\item The equation:
\[
2x+3y-6z!=222\hspace{2cm}(x,y\in\mathbb{Z},\ z\in\mathbb{N}^{0})
\]
can be solved by reducing it in modulo 3 to get $x=3k$ and by reducing
it in modulo 2 to get $y=2s$. On substituting these into the original
equation we get its (parameterized) solutions.
\end{enumerate}

\subsection{Trigonometric Equations}

The simple (and most common types) of trigonometric Diophantine equations
can be solved rather easily by investigating the cyclic behavior of
the terms of the equation and considering all the combinations (of
the values of these terms) that satisfy the given equation. The following
examples should give a reasonable insight about how to tackle and
solve (the common types) of trigonometric Diophantine equations:
\begin{enumerate}
\item The equation:
\[
3\sin\left(\frac{x\pi}{2}\right)-4\sin\left(\frac{y\pi}{2}\right)=1\hspace{2cm}(x,y\in\mathbb{Z})
\]
can be solved by noting that the sine function of integer multiples
of $\pi/2$ takes the values $0,1,0,-1$ corresponding to $x\stackrel{4}{=}0,1,2,3$
(and similarly for $y$), and hence the only possibility for the LHS
of this equation to be equal to 1 is when $\sin\left(\frac{x\pi}{2}\right)=\sin\left(\frac{y\pi}{2}\right)=-1$
so that the LHS becomes $3(-1)-4(-1)=1$. Accordingly, the solutions
of the given equation are $(x,y)=(3+4k,\ 3+4s)$ where $k,s\in\mathbb{Z}$.
\item The equation:
\[
5\sin\left(\frac{x\pi}{2}\right)+\cos(y\pi)=3\hspace{2cm}(x,y\in\mathbb{Z})
\]
can be solved by noting that the sine function of integer multiples
of $\pi/2$ takes the values $0,1,-1$ while the cosine function of
integer multiples of $\pi$ takes the values $1,-1$. Therefore, the
term $5\sin\left(\frac{x\pi}{2}\right)$ takes the values $0,5,-5$
while the term $\cos(y\pi)$ takes the values $1,-1$. As we see,
there is no combination of these values that can make the sum of the
terms on the LHS to be equal to 3. So, the given equation has no solution.
\item The equation:
\[
2\sin\left(\frac{x\pi}{2}\right)+3\cos\left(\frac{y\pi}{2}\right)=2\hspace{2cm}(x,y\in\mathbb{Z})
\]
can be solved (like the previous equation) by noting that $2\sin\left(\frac{x\pi}{2}\right)=0,2,0,-2$
corresponding to $x\stackrel{4}{=}0,1,2,3$ while $3\cos\left(\frac{y\pi}{2}\right)=3,0,-3,0$
corresponding to $y\stackrel{4}{=}0,1,2,3$. So, their sum is equal
to 2 when $(x,y)\stackrel{4}{=}(1,1)$ and $(x,y)\stackrel{4}{=}(1,3)$.
Hence, the solutions of the given equation are $(x,y)=(1+4k,\ 1+2s)$
where $k,s\in\mathbb{Z}$.
\item The equation:
\[
\sin\left(\frac{x\pi}{2}\right)+5\cos(y\pi)=6\cos(z\pi)\hspace{2cm}(x,y,z\in\mathbb{Z})
\]
can be solved by noting that:\\
\begin{tabular*}{15.95cm}{@{\extracolsep{\fill}}lll}
\noalign{\vskip0.1cm}
 & $\sin\left(\frac{x\pi}{2}\right)=0,1,0,-1$ ~ for ~ $x\stackrel{4}{=}0,1,2,3$ & \tabularnewline[0.1cm]
\noalign{\vskip0.1cm}
 & $5\cos(y\pi)=5,-5$ ~ for ~ $y\stackrel{2}{=}0,1$ & \tabularnewline[0.1cm]
\noalign{\vskip0.1cm}
 & $6\cos(z\pi)=6,-6$ ~ for ~ $z\stackrel{2}{=}0,1$ & \tabularnewline[0.1cm]
\end{tabular*}\\
So, the two sides become equal in the following two cases (which represent
the solutions of this equation):\\
\begin{tabular*}{15.95cm}{@{\extracolsep{\fill}}llll}
\noalign{\vskip0.1cm}
 & $(x,y,z)=(1+4k,\,2s,\,2t)$ & $(x,y,z)=(3+4k,\,1+2s,\,1+2t)$ & \tabularnewline[0.1cm]
\end{tabular*}\\
where $k,s,t\in\mathbb{Z}$.
\item The equation:
\[
\tan\left(x\pi+\frac{\pi}{3}\right)+2\sin\left(y\pi+\frac{2\pi}{3}\right)=2\cos\left(\frac{z\pi}{6}\right)\hspace{2cm}(x,y,z\in\mathbb{Z})
\]
can be solved by noting that:\\
\begin{tabular*}{15.95cm}{@{\extracolsep{\fill}}lll}
\noalign{\vskip0.1cm}
 & $\tan\left(x\pi+\frac{\pi}{3}\right)=\sqrt{3}$ ~~ for all $x\in\mathbb{Z}$ & \tabularnewline[0.1cm]
\noalign{\vskip0.1cm}
 & $2\sin\left(y\pi+\frac{2\pi}{3}\right)=\sqrt{3},-\sqrt{3}$ ~ for
~ $y\stackrel{2}{=}0,1$ & \tabularnewline[0.1cm]
\noalign{\vskip0.1cm}
 & $2\cos\left(\frac{z\pi}{6}\right)=2,\sqrt{3},1,0,-1,-\sqrt{3},-2,-\sqrt{3},-1,0,1,\sqrt{3}$
~ for ~ $z\stackrel{12}{=}0,1,2,\ldots,11$ & \tabularnewline[0.1cm]
\end{tabular*}\\
So, the two sides become equal in the following two cases:\\
\begin{tabular*}{15.95cm}{@{\extracolsep{\fill}}llll}
\noalign{\vskip0.1cm}
 & $(x,y,z)=(k,\,1+2s,\,3+12t)$ & $(x,y,z)=(k,\,1+2s,\,9+12t)$ & \tabularnewline[0.1cm]
\end{tabular*}\\
which can be combined in the following formula (which represents all
the solutions of the given equation): $(x,y,z)=(k,\,1+2s,\,3+6t)$
where $k,s,t\in\mathbb{Z}$.
\end{enumerate}
\clearpage{}

\section{\label{secReductionModular}Reduction and Analysis by Modular Arithmetic}

Reduction and analysis by modular arithmetic is very common method
for solving various types of Diophantine equations, and hence considering
modular arithmetic as an aiding technique for solving Diophantine
equations is highly recommended especially when dealing with polynomial
and exponential Diophantine equations. In fact, we already discussed
this method within the context of its use in solving various types
of Diophantine equations (see for instance $\S$ \ref{subModRedPoly}
and $\S$ \ref{subModRedExp}). So, in this section we discuss some
general issues about this method (namely its logical and mathematical
foundations, its purposes and some recommendations and guidelines
about its use and employment highlighting finally its importance for
non-solvability).

\subsection{\label{subLogFoundModRed}Logical and Mathematical Foundations}

Regarding the logical and mathematical foundations of the use of modular
arithmetic in solving Diophantine equations, we refer the reader to
$\S$ $\S$ 2.7.6 of \cite{SochiBookNTV1} (also see $\S$ 8.1 of
\cite{SochiBookNTV2}) where we discussed this subject\footnote{We mean the subject of relationship between the ordinary (Diophantine)
equation and the corresponding congruence equation, e.g. the relationship
between the Diophantine equation $x^{2}+3xy=0$ and the congruence
equation $x^{2}+3xy\stackrel{3}{=}0$.} in detail. So, all we need to know here are the following simple
rules:\footnote{In the following, $f(x,y)$ represents a Diophantine expression (like
$x^{2}+3xy$) in two variables noting that this applies to Diophantine
expressions in more than two variables (e.g. $x^{2}+3xy-5z^{3}$).}

\noindent $\bullet$ If $f(x,y)=0$ then $f(x,y)\stackrel{m}{=}0$
for \textit{any} $\mathbb{N}\ni m>1$, and hence (by contraposition)
if $f(x,y)\stackrel{m}{\neq}0$ for a \textit{specific} $m$ then
$f(x,y)\neq0$.

\noindent $\bullet$ If $f(x,y)\stackrel{m}{=}0$ for all $\mathbb{N}\ni m>1$
then $f(x,y)=0$.

\subsection{Purposes of Modular Reduction and Analysis}

The main (or general) purpose of using modular reduction and analysis
as an aiding tool for solving Diophantine equations is to reduce the
effort (required for searching for solution) by exploiting the versatile
collection of rules (or techniques or $\ldots$) of modular arithmetic
and its rather simplified mathematical machinery (such as the possible
use of small numbers or the classification of integers in a rather
simple and well organized residue systems). However, there are many
specific purposes that come under this general purpose. In the following
list, we try to outline some of the common specific purposes (or benefits
or objectives etc.) of using modular reduction and analysis in solving
Diophantine equations where we mostly rely for explanation and clarification
on the examples that we discussed and investigated elsewhere:
\begin{enumerate}
\item \textbf{Showing non-solvability} of the Diophantine equation (see
for instance the examples of Eqs. \ref{eqPol15x2} and \ref{eqExp4x}).
\item \textbf{Removing some variables} and hence simplifying the analysis
of the given equation (see for instance the example of Eq. \ref{eqPolx2}).
\item \textbf{Revealing patterns and clues} that can help in solving the
Diophantine equation (see for instance the examples of Eqs. \ref{eqExp5xm3y2}
and \ref{eqExp3xp5yz2}).
\end{enumerate}
We note in this regard that the employed modular reduction may not
introduce any ``visible'' change on the original Diophantine equation
apart from converting it from an ordinary equation to a congruence
equation.

\subsection{\label{subRecomModAnalysis}Recommendations and Guidelines}

We already discussed the main recommendations and guidelines about
the use of modular reduction and analysis during our discussion of
specific types of Diophantine equations (see for instance $\S$ \ref{subModRedPoly}
and $\S$ \ref{subRecomModAnalExp}). So, in the following points
we just outline and summarize what we have given earlier:

\noindent $\bullet$ Try to find and use small reduction moduli when
possible.

\noindent $\bullet$ Try to find and use a prime numbers (as a moduli)
when possible.

\noindent $\bullet$ Consider using more than one reduction modulo
in the analysis.

\noindent $\bullet$ Remember that modular reduction and analysis
is primarily a tool for showing non-solvability (see the rules of
$\S$ \ref{subLogFoundModRed} as well as $\S$ \ref{subImportanceModAnalNonSolv}).
However, it is commonly used (with further analysis and extra effort)
for investigating and finding the solutions of the Diophantine equations
(if they exist) by imposing certain restrictions and conditions on
the solutions and showing their types and forms (such as being odd
or divisible by 5).

\subsection{\label{subImportanceModAnalNonSolv}Importance of Modular Analysis
in Proving Non-Solvability}

It is important to note that modular reduction and analysis is especially
important in showing and proving that a given Diophantine equation
has no solution (either conditionally or unconditionally). In fact,
it is usually the only (or at least the main) possible tool/method
for doing this job and achieving this purpose. Therefore, as soon
as we have a guess or a hint or an indication (e.g. from initial computational
investigation; see $\S$ \ref{secInitialComputInv}) that a given
equation has no solution, we should recall our modular arithmetic
techniques and skills to find a modulo in which the given Diophantine
equation has no solution when it is reduced in that modulo. In fact,
parity analysis and check (which is one of the most common and basic
methods for establishing non-solvability; see for instance \ref{secParityChecks})
is no more than modular reduction and analysis in modulo 2 (where
its wide spread use and its distinction from modular analysis is because
of its simplicity and intuitivity as well as other reasons and factors).

\clearpage{}

\section{\label{secFactorization}Factorization Analysis}

Factorization analysis is a very common and useful method for analyzing
and solving Diophantine equations and hence it should always be considered
as one of the first approaches in tackling Diophantine problems. In
its most common and simple form (noting that it has various forms
and variants as will be outlined next) factorization analysis is based
on our ability to produce a factored expression involving variables
(on the LHS) that is equal to a specific number (on the RHS). The
number of the possible factors on the two sides may also be considered
in some variants of factorization analysis since it can eliminate
certain possibilities for the solution. Other factors and considerations
(such as divisibility or primality/composity) may also be considered
and included in more versatile factorization analysis arguments.

In the following subsections we will briefly investigate the aforementioned
variants of factorization analysis.

\subsection{Simple Factorization Analysis}

As indicated in the preamble of this section, simple factorization
analysis is based on producing a factored expression involving variables
(on the LHS) that is equal to a specific number (on the RHS) where
the factors involving variables (on the LHS) can be matched with numeric
factors (on the RHS) to produce systems of simultaneous equations
that can be solved to produce the solution(s) of the given problem.
In fact, some examples of this type of factorization analysis have
already been given (see for instance the examples of point \ref{enuFCp1}
of $\S$ \ref{subFactorizationComparison}).

\subsection{Factorization Analysis based on the Number of Factors}

Factorization analysis may also be based on comparing the number of
factors (i.e. not their specific form or value) on the two sides of
the equation. For example, the function:
\[
f(x,y)=120x^{5}+274x^{4}y+225x^{3}y^{2}+85x^{2}y^{3}+15xy^{4}+y^{5}
\]
can be factorized in the following 5-factor form:
\[
f(x,y)=(x+y)(2x+y)(3x+y)(4x+y)(5x+y)
\]
and hence we can conclude (through factorization analysis based on
the number of factors on the two sides) that $f(x,y)=21$ has no solution
(because 21 can be factorized only in 2 or 3 or 4 distinct integer
factors) while $f(x,y)=45$ can have (and actually has) solution (because
45 can be factorized into 5 distinct integer factors). However, it
should be noticed that this type of factorization analysis is primarily
for proving non-solvability and hence solvability (and obtaining the
solutions) requires extra work (noting that having the same number
of factors on the two sides of equation is a necessary but not sufficient
condition for solvability in this type of factorization analysis).

\subsection{Factorization with Divisibility Analysis}

In this type of factorization analysis we consider employing divisibility
arguments of one side (or of factors of one side) by the other side
(or by factors of the other side). In fact, some examples of this
type of factorization analysis have already been given (see for instance
the examples of points \ref{enuFCp2} and \ref{enuFCp3} of $\S$
\ref{subFactorizationComparison}).

\subsection{Factorization with Primality/Composity Analysis}

In this type of factorization analysis we consider employing primality/composity
arguments where one side (or factors of one side) is compared to the
other side (or to factors of the other side) from this perspective
(i.e. being prime or composite). Some simple examples of this type
of factorization analysis have been given in $\S$ \ref{subPrimalityComposityChecks}.

\clearpage{}

\section{\label{secComparisonToSimilar}Comparison to Similar Equations and
Problems}

An ideal way for solving a given Diophantine equation is to compare
it to a similar equation whose solutions are known (or whose solutions
are easier to obtain) and hence obtain the solutions of the given
equation with no effort (or with minimal effort). In fact, sometimes
solving a given Diophantine equation may require no more than copying
and pasting the solution of a previously-solved problem with some
modifications and adaptations to reflect the specific characteristics
of the given Diophantine equation. Therefore, it is recommended when
tackling a Diophantine equation to consider this method of solution
as one of the first options to recall and consider.

As hinted above, the advantage of the method of comparison is simplicity
and ease where the effort required to solve the given Diophantine
equation is reduced substantially. However, an obvious limitation
of this method is that we need to have a similar equation with known
(or easy-to-obtain) solutions which is obviously not available in
most cases. However, it is useful to search for such an equation in
our collection of previously-solved Diophantine equations, and hence
keeping an organized and classified ``database'' of Diophantine
equations (with known solutions) is very useful (especially to those
who specialize in the subject of Diophantine equations).\footnote{In fact, this is one reason for organizing and classifying Diophantine
equations in our books (see \cite{SochiBookNTV1,SochiBookNTV2}).} Some examples of this method are given in the following:
\begin{enumerate}
\item Let us assume that we already investigated the Diophantine equation:
\begin{equation}
5x^{3}+4y^{3}-9=0\label{eqasdsd1}
\end{equation}
and obtained its solutions which are $(x,y)=(1,1)$ and $(x,y)=(13,-14)$.
Now, if we have to solve the Diophantine equation:
\begin{equation}
5x^{3}-4y^{3}+9=0\label{eqasdsd2}
\end{equation}
then all we need to do is to change the sign of the $x$ value of
the solutions of Eq. \ref{eqasdsd1} and hence obtain the solutions
$(x,y)=(-1,1)$ and $(x,y)=(-13,-14)$. This is because if we multiply
Eq. \ref{eqasdsd2} by $-1$ then we obtain:
\[
5(-x)^{3}+4y^{3}-9=0
\]
which is no more than Eq. \ref{eqasdsd1} with the sign of $x$ being
reversed.\\
Similarly, if we have to solve the Diophantine equation:
\begin{equation}
5x^{3}+4y^{3}+9=0\label{eqasdsd3}
\end{equation}
then all we need to do is to change the sign of the $x$ and $y$
values of the solutions of Eq. \ref{eqasdsd1} and hence obtain the
solutions $(x,y)=(-1,-1)$ and $(x,y)=(-13,14)$. This is because
if we multiply Eq. \ref{eqasdsd3} by $-1$ then we obtain:
\[
5(-x)^{3}+4(-y)^{3}-9=0
\]
which is no more than Eq. \ref{eqasdsd1} with the sign of $x$ and
$y$ being reversed.
\item Let us assume that we already investigated the Diophantine equation:
\begin{equation}
x^{3}+y^{3}+x^{2}+y^{2}=0\label{eqrtrhh1}
\end{equation}
and obtained its solutions which are $(x,y)=(0,0)$, $(0,-1)$, $(-1,0)$,
$(-1,-1)$. Now, if we have to solve the Diophantine equation:
\begin{equation}
x^{3}+y^{3}-x^{2}-y^{2}=0\label{eqrtrhh2}
\end{equation}
then all we need to do is to change the sign of the $x$ and $y$
values of the solutions of Eq. \ref{eqrtrhh1} and hence obtain the
solutions $(x,y)=(0,0)$, $(0,1)$, $(1,0)$, $(1,1)$. This is because
if we multiply Eq. \ref{eqrtrhh2} by $-1$ then we obtain:
\[
(-x)^{3}+(-y)^{3}+x^{2}+y^{2}=0\hspace{2cm}\textrm{i.e.}\hspace{2cm}(-x)^{3}+(-y)^{3}+(-x)^{2}+(-y)^{2}=0
\]
which is no more than Eq. \ref{eqrtrhh1} with the sign of $x$ and
$y$ being reversed.
\item We should also refer to the examples of $\S$ \ref{subComparison}
and $\S$ \ref{subCompareFrac} where we used the method of comparison
for solving some Diophantine equations of certain types.
\end{enumerate}
We should finally note that comparison to similar equations and problems
may require some manipulations and transformations (see the examples
in $\S$ \ref{secManipulationsTransformation}). In fact, we may even
extend the method of comparison beyond direct comparison of two similar
equations with a specific form and hence we may consider comparing
equations of certain characteristic features (though they may not
look similar in form) that make their method of solution (or the rationale
and logic behind their method of solution) similar.

\clearpage{}

\section{\label{secSymmetriesCycling}Symmetries and Cycling Patterns in the
Variables}

Symmetry in the variables\footnote{``Symmetry'' in the variables of an algebraic expression means that
the variables can be exchanged without affecting the expression.} of a Diophantine equation is a useful feature that can be exploited
in assuming temporarily that the variables have a certain increasing
or decreasing order. This should facilitate the search for a solution
where the final and complete solution can be obtained eventually by
permuting the solution obtained on the base of ordering assumption.
So, when tackling a Diophantine equation with total or partial symmetry
in its variables it is recommended to consider such a symmetry by
synthesizing a logical or mathematical argument that exploits such
a symmetry.

We should also look for any cyclic pattern\footnote{``Cyclic pattern'' in the variables of an algebraic expression means
that the variables can be exchanged in a certain cyclic order without
affecting the expression although the expression is affected if the
variables are exchanged without regard to the cyclic order. So, cyclic
pattern is a restricted form of symmetry.} in the variables of Diophantine equations where this pattern can
be exploited similarly, i.e. by assuming initially that a certain
variable is the biggest or the smallest (since cycling would allow
us to bring this variable to a certain position in the equation where
ordering can be exploited to make an argument that leads to the interim
solution and this solution can be generalized later by lifting the
condition of ordering).

The following are some examples of how symmetry and cyclic pattern
in the variables can be exploited in solving Diophantine equations:
\begin{enumerate}
\item The Diophantine equation:
\[
\frac{1}{x}+\frac{1}{y}=\frac{2}{3}\hspace{2cm}(x,y\in\mathbb{N})
\]
is symmetric in $x,y$ and hence we can assume initially that $x\leq y$.
Accordingly, we can argue (based on this assumption and considering
the magnitude of the RHS as well as similar factors) that $x$ must
be either 2 (and hence $y=6$) or 3 (and hence $y=3$). So, we obtain
the interim solutions $(x,y)=(2,6)$ and $(3,3)$. Now, if we lift
the condition $x\leq y$ (thanks to the symmetry) we obtain (by permutation)
another solution, i.e. $(x,y)=(6,2)$. So, the given Diophantine equation
has only three solutions.
\item The Diophantine equation:
\[
\frac{1}{x}+\frac{1}{y}+\frac{1}{z}=2\hspace{2cm}(x,y,z\in\mathbb{N})
\]
is symmetric in $x,y,z$ and hence we can assume initially that $x\leq y\leq z$.
Accordingly, we can argue (based on this assumption) that $x$ cannot
be greater than 1 and $y$ must be 2 and hence $z$ must be 2. This
argument produces the interim solution $(x,y,z)=(1,2,2)$ which can
be generalized by permutation (thanks to the symmetry) to produce
the other two solutions, i.e. $(x,y,z)=(2,1,2)$ and $(2,2,1)$.
\item The Diophantine equation:
\[
x+y+z=xyz\hspace{2cm}(x,y,z\in\mathbb{Z},\ xyz\neq0)
\]
is symmetric in $x,y,z$ and hence we can assume initially that $|x|\leq|y|\leq|z|$.
Accordingly, we can argue (based on this assumption) that $x$ must
be equal to $\pm1$ and hence either $y=-2$ and $z=-3$ (i.e. when
$x=-1$) or $y=2$ and $z=3$ (i.e. when $x=1$). This argument produces
the two interim solutions $(x,y,z)=\pm(1,2,3)$ which can be generalized
by permutation (thanks to the symmetry) to produce the other ten solutions,
i.e. $(x,y,z)=\pm(1,3,2)$, $\pm(2,1,3)$, $\pm(2,3,1)$, $\pm(3,1,2)$,
$\pm(3,2,1)$.
\item The Diophantine equation:
\[
x^{2}+y^{2}=z^{2}\hspace{2cm}(x,y,z\textrm{ are consecutive natural numbers})
\]
is symmetric in $x,y$ and hence we can assume initially that $x<y$
(i.e. $y=x+1$), that is:\\
\begin{tabular*}{15.95cm}{@{\extracolsep{\fill}}lllll}
\noalign{\vskip0.1cm}
$x^{2}+(x+1)^{2}=(x+2)^{2}$ & $\rightarrow$ & $x^{2}-2x-3=0$ & $\rightarrow$ & $(x+1)(x-3)=0$\tabularnewline[0.1cm]
\end{tabular*}\\
i.e. $x=3$ (noting that $x\in\mathbb{N}$). Therefore, we have only
one interim solution, i.e. $(x,y,z)=(3,4,5)$. Now, if we lift the
condition $x<y$ (thanks to the symmetry) we obtain (by permutation)
the other solution, i.e. $(x,y,z)=(4,3,5)$.
\item The Diophantine equation:
\[
x+y+z+w=xyzw\hspace{2cm}(x,y,z,w\in\mathbb{N})
\]
is symmetric in $x,y,z,w$ and hence we can assume initially that
$x\leq y\leq z\leq w$. Accordingly, we can argue as before (based
on this assumption) that we must have $x=y=1$ and $z=2$ and hence
$w=4$. So, the interim solution is $(x,y,z,w)=(1,1,2,4)$ which can
be generalized by permutation (thanks to the symmetry) to produce
the other eleven solutions.
\item The Diophantine equation $x^{4}+y^{4}+z^{4}=3042$ (which will be
discussed in $\S$ \ref{secUpperLower}) is another example of symmetry
in the variables.
\item The Diophantine equation:
\[
x^{3}y+y^{3}z+z^{3}x=1\hspace{2cm}(x,y,z\in\mathbb{N}^{0})
\]
is not symmetric in its variables, e.g. if we exchange $x$ and $y$
we get:
\[
y^{3}x+x^{3}z+z^{3}y=1
\]
which is not the same as the original equation. However, it is cyclic
in its variables, i.e. if we exchange $x\rightarrow y$, $y\rightarrow z$
and $z\rightarrow x$ we get:
\[
y^{3}z+z^{3}x+x^{3}y=1
\]
which is the same as the original equation (apart from the order of
terms which is irrelevant in this context). So, we can assume initially
that $x\leq y$ and $x\leq z$ (since we can cycle the variables to
put $x$ in such a position in the equation noting that we cannot
assume that $x\leq y\leq z$ because such an ordering flexibility
requires full symmetry). Accordingly, we can argue (based on this
assumption) that $x$ must be 0 (because otherwise the LHS will be
greater than 1) and hence we must have $y=z=1$. So, the interim solution
is $(x,y,z)=(0,1,1)$. Now, if we cycle the variables we obtain the
other two solutions, i.e. $(x,y,z)=(1,0,1)$ and $(x,y,z)=(1,1,0)$.
\end{enumerate}
\clearpage{}

\section{\label{secUpperLower}Upper and Lower Bounds}

Imposing upper or/and lower bounds or limits (usually in the form
of bounding inequalities) on the potential solutions should be considered
when tackling Diophantine equations (e.g. when dealing with equations
involving fractions). This can reduce the complexity of the process
of search for solution substantially since it reduces the number of
possible solutions. In fact, it can reduce this number from being
infinite to finite and this should allow the use of more simple methods
of inspection and search (e.g. by computational tools and techniques).

The two main considerations for imposing such bounds are:

\noindent $\bullet$ Magnitude\footnote{``Magnitude'' here should mean the position on the number line rather
than the absolute value (as we usually use).} considerations where the potential solutions cannot exceed in magnitude
certain (upper or/and lower) limits. Bounds imposed by magnitude considerations
can impose a limit on the number of possible solutions, i.e. they
can make this number finite when both lower and upper bounds are imposed.

\noindent $\bullet$ Sign considerations where the potential solutions
can be of only one type of sign (i.e. positive or negative). Although
bounds imposed by sign considerations should reduce the difficulty
of the problem they usually do not impose a limit on the number of
possible solutions, i.e. this number remains infinite. In fact, sign
considerations can be seen as a special form of magnitude considerations
(which we discussed in the previous point) where only upper/lower
bound is imposed.

The following are a few examples for the use of bounds (and bounding
arguments) in solving Diophantine problems:
\begin{enumerate}
\item Consider the following Diophantine equation:
\[
\frac{1}{x}+\frac{1}{y}=z\hspace{2cm}(x,y,z\in\mathbb{Z},\ xy\neq0)
\]
It should be obvious that the LHS (and hence $z$) cannot be less
than $-2$ or greater than $+2$ and hence:\\
$\bullet$ If $z=-2$ then $(x,y,z)=(-1,-1,-2)$.\\
$\bullet$ If $z=-1$ then $(x,y,z)=(-2,-2,-1)$.\\
$\bullet$ If $z=0$ then $(x,y,z)=(k,-k,0)$ where $\mathbb{Z}\ni k\neq0$.\\
$\bullet$ If $z=1$ then $(x,y,z)=(2,2,1)$.\\
$\bullet$ If $z=2$ then $(x,y,z)=(1,1,2)$.
\item Consider the following Diophantine equation:
\[
x^{4}+y^{4}+z^{4}=3042\hspace{2cm}(x,y,z\in\mathbb{Z})
\]
It should be obvious that any one of the three variables cannot be
less than $-7$ or greater than $+7$ (because otherwise the LHS will
be greater than 3042), and hence if we carry a simple computational
search within the limits $-7\leq x,y,z\leq+7$ we will find all the
(forty eight) solutions of this equation.
\item The following Diophantine equation:
\[
x^{2}-xy+y^{2}+2x-y=2\hspace{2cm}(x,y\in\mathbb{Z})
\]
can be put in the following form:
\[
(x+2)^{2}+(y-1)^{2}+(x-y)^{2}=9
\]
and hence we must have $-3\leq(x+2)\leq3$, i.e. $-5\leq x\leq1$.
Now, if we investigate all the seven possibilities (i.e. $x=-5,-4,-3,-2,-1,0,1$)
we get all the (six) solutions of the given equation.
\item Consider the following Diophantine equations:
\[
2^{x}+3^{y}=1\hspace{3cm}4^{x}+9^{y}=2\hspace{2cm}(x,y\in\mathbb{N}^{0})
\]
It should be obvious that the first equation has no solution because
its LHS cannot be less than 2, while the second equation has only
the trivial solution $x=y=0$ because otherwise the LHS will be greater
than 2.
\item Consider the following Diophantine equations:
\[
\frac{1}{x}+\frac{2}{y}+\frac{3}{z}=7\hspace{3cm}\frac{x}{y}+\frac{y}{z}+\frac{z}{x}=2\hspace{2cm}(x,y,z\in\mathbb{N})
\]
It should be obvious that these equations have no solution because
the LHS of the first equation cannot exceed 6 while the LHS of the
second equation must exceed 2.
\item Consider the following Diophantine equation:
\[
z=\frac{1}{x}+\frac{1}{y}+\frac{2}{xy}\hspace{2cm}(x,y,z\in\mathbb{N})
\]
It should be obvious that the RHS (and hence $z$) cannot exceed 4.
By a similar magnitude argument it can be easily shown that $x$ and
$y$ cannot exceed 4. So, all we need to do is to inspect (computationally)
the combinations of $x,y,z=1,2,3,4$ to obtain the (five) solutions
of the given equation.
\item Consider the following system of Diophantine equations:
\begin{equation}
3xy+5x^{3}y=230\hspace{1.5cm}\textrm{and}\hspace{1.5cm}x^{2}+xy=14\hspace{2cm}(x,y\in\mathbb{Z})\label{eqDioSys1}
\end{equation}
From the first equation it is fairly obvious that $x$ and $y$ must
have the same sign (because otherwise the LHS will be non-positive).
Now, from the second equation it should be obvious that the absolute
value of $x$ cannot exceed 3 and hence we have only the following
6 possibilities to consider: $x=-3,-2,-1,1,2,3$. On inspecting these
possibilities we find the solutions of this system, i.e. $(x,y)=(-2,-5)$
and $(x,y)=(2,5)$.
\end{enumerate}
\clearpage{}

\section{Reduction of Domain}

When dealing with Diophantine equations it is recommended to consider
reducing the domain of solution temporarily until a solution is found
(for the reduced domain) where this solution can be extended and generalized
later on to reach the final and complete solution for the entire domain.
For example, if we are dealing with a Diophantine equation in the
domain of integers $\mathbb{Z}$ involving variables with even powers
then we can start by considering its solution in the domain of natural
numbers $\mathbb{N}$ (instead of $\mathbb{Z}$) where the final and
complete solution can be obtained later on by extending the domain
to the negative integers (noting that even powers do not distinguish
between positive and negative bases).

The reduction of domain may also take the form of splitting the domain
into parts and dealing with these parts separately (and usually independently).
For example, we may consider splitting\footnote{Such a split can involve some or all variables (depending on the nature
of the Diophantine problem and the tackling strategy).} the integer domain of a Diophantine equation to negative integers
and positive integers where we investigate these sub-domains separately
(e.g. by applying different arguments and techniques to each sub-domain)
and the final solution (over the entire integer domain) will be obtained
in the end by taking the union of the solutions in the sub-domains
(noting the possibility of having no solution in some or all sub-domains).
Such a scenario may also be considered by splitting the domain to
odd integers and even integers (or to integers greater than and integers
less than a certain value or $\ldots$) where we deal with these sub-domains
separately and obtain the entire solution in stages (as before).

Anyway, reduction/splitting of domain in its various shapes and forms
is a very common approach in tackling and solving Diophantine problems
and hence it should always be considered since it usually brings many
benefits (such as reducing the complexity of the given problem and
making the strategy of tackling it tidy and organized) and can provide
key clues for solving it. It can also provide a partial solution to
the given problem, i.e. when the investigation leads to solution only
in some parts of the domain (where this partial solution can be extended
in the future to include the entire domain or where the original problem
is modified to consider such a restriction on the domain).

We present in the following some examples for the reduction of domain
in tackling Diophantine problems:
\begin{enumerate}
\item The ``Pythagorean equation'' $x^{2a}+y^{2b}=z^{2c}$ (where $x,y,z\in\mathbb{Z}$
and $a,b,c\in\mathbb{N}$) is an obvious example for the possibility
of considering the reduction of domain where we can start by considering
initially the solutions in the reduced domain of natural numbers (or
non-negative numbers). Any solution obtained in this restricted domain
can then be extended to the domain of integers by considering all
the possible sign alterations in the values of the variables in the
natural solutions or non-negative solutions).
\item Another obvious example for the reduction of domain is the Diophantine
equations of the Pell type (see $\S$ \ref{subPellEquation}) where
we can obtain the solutions in the reduced domain of natural numbers
first (e.g. by the use of standard Pell's techniques) and extend these
solutions subsequently to the domain of integers by considering all
the possible sign alterations in the values of the variables in the
natural solutions.
\item Consider the following Diophantine equation:
\[
x^{3}y+3xy^{3}+7xy=1085\hspace{2cm}(x,y\in\mathbb{Z})
\]
It is fairly obvious that $x$ and $y$ have the same sign (because
otherwise the LHS will be non-positive). Now, if we start by assuming
that $x,y\in\mathbb{N}$ then we can easily obtain the solution $(x,y)=(1,7)$.
Following this, we extend this reduced domain (by lifting the condition
$x,y\in\mathbb{N}$) to obtain the other solution, i.e. $(x,y)=(-1,-7)$
noting that the equation does not change by reversing the sign of
$x$ and $y$.
\item Consider the following Diophantine equation:
\[
\frac{x}{y}+xy=2\hspace{2cm}(x,y\in\mathbb{Z},\ y\neq0)
\]
It is obvious that $x$ and $y$ have the same sign (because otherwise
the LHS will be non-positive). Now, if we start by assuming that $x,y\in\mathbb{N}$
then we can easily obtain the solution $(x,y)=(1,1)$. Following this,
we extend this reduced domain (by lifting the condition $x,y\in\mathbb{N}$)
to obtain the other solution, i.e. $(x,y)=(-1,-1)$ noting that the
equation does not change by reversing the sign of $x$ and $y$.
\item Consider the system of Eq. \ref{eqDioSys1}. Noting that $x$ and
$y$ must have the same sign, we can start by assuming that $x,y\in\mathbb{N}$
where we can easily obtain the solution $(x,y)=(2,5)$. We then extend
this reduced domain (by lifting the condition $x,y\in\mathbb{N}$)
to obtain the other solution, i.e. $(x,y)=(-2,-5)$.\footnote{Noting that $x$ and $y$ have the same sign (as well as $x^{2}$
in the second equation), it should be obvious that the equations of
this system do not change by such a change in sign and hence the second
(i.e. negative) solution can be obtained by just reversing the sign
of the first (i.e. positive) solution.}
\end{enumerate}
\clearpage{}

\section{\label{secBasicRules}Basic Rules and Principles}

When we deal with a Diophantine problem we should always try to recall
common rules and basic principles that can help in solving the problem
and reducing the amount of work required to solve it (since such rules
and principles usually summarize and encompass a number of steps in
the formal proof or argument required for solving the problem). In
fact, we should even consider manipulating the given Diophantine equations
in such a way that enables us to exploit such rules and principles.
Prominent examples of such rules and principles are the (integer)
ordering rules, the algebraic rules of (integer) powers, and the rules
of primality and coprimality. These rules and principles are particularly
useful for identifying the Diophantine equations and systems that
have no solution (or have restrictions on their solutions).

A few examples of these rules and principles are given in the following
points (noting that there are many other rules and principles like
these that can be exploited for tackling Diophantine problems):

\noindent $\bullet$ No (perfect) square can be between consecutive
(perfect) squares, no (perfect) cube can be between consecutive (perfect)
cubes, and so on.

\noindent $\bullet$ No (perfect) square can be the sum of two odd
squares.\footnote{This rule is actually based on the rules and properties of Pythagorean
triples (see $\S$ \ref{subPythagoreanTriples}) which will be mentioned
later in this list.}

\noindent $\bullet$ Any (perfect) cube can be the difference of two
(perfect) squares.

\noindent $\bullet$ The difference between two consecutive cubes
cannot be divisible by 5.

\noindent $\bullet$ Any integer can be the sum of two perfect squares
minus a perfect square.

\noindent $\bullet$ Any natural number can be the sum of four non-negative
integer squares.

\noindent $\bullet$ The difference between two non-trivial perfect
squares cannot be 4.

\noindent $\bullet$ If the square root of an integer is rational
then the integer is a perfect square (i.e. the square root is an integer).

\noindent $\bullet$ Any two consecutive integers are coprime.

\noindent $\bullet$ The natural powers of distinct primes are coprime.

\noindent $\bullet$ Any odd prime is congruent (in modulo 4) either
to $+1$ or to $-1$ (or similarly to $+3$).

\noindent $\bullet$ In any three consecutive odd integers exactly
one of them is divisible by 3.

\noindent $\bullet$ All factorials are even except $0!$ and $1!$.

\noindent $\bullet$ The factorial $n!$ is not divisible by any prime
$p>n$.

\noindent $\bullet$ All binomial coefficients and multinomial coefficients
are integers (despite their appearance as ratios or fractions).

\noindent $\bullet$ The rules (and properties) of Pythagorean triples
should also be considered in this regard.

\noindent $\bullet$ Fermat's last theorem (see $\S$ \ref{subFermatLastTheorem})
can also be considered in this regard (since it can be seen as a rule
that eliminates the possibility of existence of solution or the possibility
of existence of certain types of solution).

We present in the following some examples for the use of such rules
and principles in tackling and solving Diophantine problems:
\begin{enumerate}
\item By simple algebraic manipulation, the following Diophantine equation:
\[
4x^{2}+16y^{2}-9z^{2}-12x+8y+10=0\hspace{2cm}(x,y,z\in\mathbb{Z})
\]
can be put in the following form:
\[
(2x-3)^{2}+(4y+1)^{2}=(3z)^{2}
\]
Now, if we recall the rule that ``the sum of two odd squares cannot
be a perfect square'' then we can conclude that this Diophantine
equation has no solution (i.e. the problem is solved with virtually
no effort thanks to this rule).
\item The following Diophantine equation:
\[
216x^{3}+27y^{3}+216x^{2}+72x-721=0\hspace{2cm}(x,y\in\mathbb{N})
\]
can be put in the following form:
\[
(6x+2)^{3}+(3y)^{3}=9^{3}
\]
and hence it can be solved by Fermat's last theorem, i.e. it has no
solution in $\mathbb{N}$ (and in fact it has no solution even in
$\mathbb{Z}$).
\item The following Diophantine equation:
\[
16x^{4}+81y^{4}=z^{4}\hspace{2cm}(x,y,z\in\mathbb{Z})
\]
can be written as:
\[
(2x)^{4}+(3y)^{4}=z^{4}
\]
and hence by Fermat's last theorem we can conclude that it can have
only trivial solutions (i.e. solutions with $xyz=0$). This restriction
on the solutions should reduce the difficulty of searching for solutions
substantially.
\item Consider the following system of Diophantine equations:
\[
y^{3}+3y^{2}-x^{2}+3y-z^{5}+z+1=0\hspace{1.5cm}\textrm{and}\hspace{1.5cm}x^{2}-y^{3}=0\hspace{2cm}(x,y,z\in\mathbb{Z})
\]
If we put the equations of this system in the following forms:
\[
(y+1)^{3}-x^{2}=z^{5}-z\hspace{1.5cm}\textrm{and}\hspace{1.5cm}x^{2}=y^{3}
\]
and substitute from the second equation into the first equation then
we get:
\[
(y+1)^{3}-y^{3}=z^{5}-z
\]
Now, if we remember the rule that ``the difference between two consecutive
cubes cannot be divisible by 5'' and note that $z^{5}-z$ is divisible
by 5 then we can conclude that this system has no solution.
\end{enumerate}
\clearpage{}

\section{\label{secOtherRecommendations}Other Recommendations}

There are many other recommendations that should be considered when
dealing with Diophantine problems (noting that many of these recommendations
usually depend on the particular type of Diophantine equation/system
of concern). Examples of these recommendations are:
\begin{enumerate}
\item Recall interesting theorems (such as Wilson's theorem and Fermat's
little theorem) possibly as part of modular arithmetic investigation
and analysis to the Diophantine equation. For example, Wilson's theorem
may be useful to recall and use (within modularity investigation and
analysis) when the Diophantine problem involves factorials.
\item Consider special (or limiting or obvious or eccentric or $\ldots$)
cases and instances such as when one (or more) of the variables is
0 or $\pm1$ or goes to infinity or becomes negative. Such considerations
can give an insight in the solution (or reduce the possibilities or
organize the approach of solution or give a clue to the solution or
$\ldots$).
\item Give special attention to the dominating terms in the equation(s)
which can (for instance) impose limits or determine the eventual tendency
of the equation(s).
\item Consider reformulating the problem that is at the base of the given
Diophantine equation/system such that a solution can be obtained.
For instance, we may impose certain limits and restrictions on the
domain of the problem which enable us to obtain full solution.
\item Consider obtaining partial (or tentative or conditional or $\ldots$)
solutions when full (or certain or unconditional or $\ldots$) solutions
could not be reached. Such ``interim'' solutions can form a basis
for future attempts and investigations to obtain final solutions.\footnote{This recommendation (and the previous one) takes into consideration
practical factors and issues. In brief, we should not give up! Moreover,
we should always try to make use of what we could obtain and achieve
(e.g. partial or conditional solution). This helps to maintain confidence
and high spirit which help us in our future investigations of Diophantine
problems. Such psychological factors are very important for our success
in solving Diophantine problems. In fact, Diophantine equations is
not an easy subject and hence without such confidence and high spirit
our ability to solve these equations will be reduced and diminished.}
\end{enumerate}
\clearpage{}

\section{Recommendations for Systems of Equations}

Systems of Diophantine equations can be classified into two main categories:
linear (when all the equations of the system are linear) and non-linear
(when some or all of the equations of the system are non-linear).
These two categories will be investigated in the following two subsections.
However, before we go through this investigation we would like to
draw the attention to the following points:
\begin{enumerate}
\item Many of the (previously discussed) recommendations related to individual
Diophantine equations applies equally or similarly to systems of Diophantine
equations as well (at least by considering the individual equations
as part of the system but not the system itself). In fact, we already
indicated (frequently) such system recommendations (at least within
the given examples for certain recommendations and guidelines).
\item We should remind the reader of what we have said before (see $\S$
\ref{secIntroduction}) that is: there are two main methods for solving
systems of Diophantine equations in number theory. The first is based
on using the traditional methods of solving systems of multivariate
equations (as investigated in algebra and linear algebra for instance)
such as by substitution or comparison or use of the techniques of
matrices, and the second is by solving the individual equations separately
(either by the general methods of algebra or by the special methods
and techniques of number theory) and selecting the solutions that
satisfy the system as a whole (i.e. by accepting only the solutions
which are common to all the equations).
\end{enumerate}

\subsection{Systems of Linear Equations}

Systems of linear Diophantine equations can be easily solved by the
well known (and standard) methods of linear algebra where only the
integer (or sub-integer such as positive integer) solutions are accepted.
Other methods (including the methods used for solving systems of non-linear
equations) are also possible to use in general for solving systems
of linear equations (as will be discussed later).

\subsection{Systems of Non-Linear Equations}

Non-linear systems\footnote{``Non-linear systems'' should be more accurate than ``systems of
non-linear equations'' if we consider systems that consist of some
linear and some non-linear equations (which should be classified as
non-linear systems although some of their equations are linear).} of Diophantine equations can be solved usually by the well known
algebraic methods of solving systems of equations (which are not restricted
to Diophantine equations) such as by substitution, elimination and
comparison (see $\S$ \ref{subSubstitution}) where only integer (or
sub-integer) solutions are accepted. However, it should be noted that
certain methods (which are mostly investigated in linear algebra)
are not applicable in solving non-linear Diophantine systems since
these methods are specific to linear systems. In the following sub-subsections
we will investigate briefly some of the common recommendations and
guidelines (as well as methods and techniques) for solving non-linear
systems of Diophantine equations. However, before that we would like
to draw the attention to the following useful remarks:
\begin{enumerate}
\item It is recommended to test all the obtained (integer) solutions on
the original system of equations. This is because some of the algebraic
manipulations which are required during the process of solving systems
of non-linear equations (such as raising to powers or multiplication
or division by a variable) can introduce foreign solutions and hence
by testing the obtained solutions we make sure that no foreign solution
is introduced during these manipulations.
\item Systems of linear equations can be considered as a special type of
systems of non-linear equations and hence they can be solved by the
methods of systems of non-linear equations (when applicable) as well
as by the methods and techniques which are specific to systems of
linear equations (such as by the methods of matrices of linear algebra).
\item Unlike systems of linear Diophantine equations, there is no standard
(or systematically applicable) method or technique for solving systems
of non-linear Diophantine equations and hence solving these systems
is a mix of art and science (where the art is required for selecting
the main approach or strategy for solving the system while the science
is needed for employing and applying the technicalities required by
the selected approach).\footnote{In fact, even solving individual Diophantine equations is mostly a
mix of art and science.} Therefore, it is recommended to investigate various potential methods
or strategies before setting off and making the choice of the best
(or even applicable) strategy to use. Investing some time and effort
in this initial investigation can be very rewarding and beneficial
and can save considerable amount of time and effort in trying to solve
the system in a rather random approach and chaotic manner.
\end{enumerate}

\subsubsection{Initial Sensibility Checks}

As indicated in $\S$ \ref{subFinThoughtInitSensCheck}, initial sensibility
checks are recommended (and even required) as the first step in tackling
and solving Diophantine systems as well as in tackling and solving
individual Diophantine equations. In many cases, solving the given
system of equations does not need more than these initial sensibility
checks since these checks can reveal that the system either has no
solution (e.g. because one of the equations has no solution due to
parity violation or modularity inconsistency or because some of the
equations require conditions that contradict the conditions required
by the other equations) or because the solution becomes so obvious
by these initial sensibility checks.

For example, initial sensibility checks should reveal that the following
system (where $x,y\in\mathbb{Z}$):

\begin{tabular*}{15.95cm}{@{\extracolsep{\fill}}llllll}
\noalign{\vskip0.1cm}
$3xy+x^{2}-5y=33$ & and & $5x^{3}+10xy+11y^{2}=23$ & and & $17x+4xy^{2}=147$ & \tabularnewline[0.1cm]
\end{tabular*}

\noindent has no solution because the second equation has no solution
(noting that in modulo 5 this equation becomes $y^{2}\stackrel{5}{=}3$
which has no solution since 3 is not a quadratic residue of 5).

On the other hand, initial sensibility checks (or inspection) should
reveal that the following system:

\begin{tabular*}{15.95cm}{@{\extracolsep{\fill}}llllll}
\noalign{\vskip0.1cm}
 & $7^{x}-8^{y}=48$ & and & $3x^{2}-5y^{3}=12$ & $(x,y\in\mathbb{N}^{0})$ & \tabularnewline[0.1cm]
\end{tabular*}

\noindent has (only) the obvious solution $(x,y)=(2,0)$. This is
because according to the first equation $y$ must be 0 (to avoid parity
violation) and hence $(x,y)=(2,0)$ is the only possible solution
to this system.

Also see the examples given in $\S$ \ref{subFinThoughtInitSensCheck}.

\subsubsection{Graphic Investigation and Reasoning}

Some systems can be easily solved by graphic investigation and reasoning
(which may or may not require plotting of actual graphs representing
the equations of the system). For example, the following system of
Diophantine equations:

\begin{tabular*}{15.95cm}{@{\extracolsep{\fill}}llllll}
\noalign{\vskip0.1cm}
 & $x^{2}-2x+4+y=0$ & and & $x^{2}+y^{2}+6x-10y+30=0$ & ($x,y\in\mathbb{Z}$) & \tabularnewline[0.1cm]
\end{tabular*}

\noindent was solved in $\S$ \ref{secPrelimGraphInspect} by graphical
reasoning without need for plotting any graph. Many other systems
can be similarly solved either graphically (by plotting actual graphs)
or by pure graphical reasoning (without plotting any graph). So, it
is recommended to consider graphic investigation when tackling Diophantine
systems (especially non-linear systems with two variables).

\subsubsection{Test of Solutions of Known Equation}

If the system of Diophantine equations contains an equation whose
solution is known (or can be obtained easily or more easily) then
the best approach for solving the system is to test the solutions
of that equation on the other equations where only the common solutions
(if any) to all equations are accepted. This is based on the obvious
fact that the solution of the system is the intersection of the solutions
of the individual equations (see the paragraph before the last of
$\S$ \ref{secIntroduction}) and hence the set of solutions of the
system cannot exceed the set of solutions of any one of the equations
in the system, i.e. the set of solutions of the system is a (proper
or improper) subset of the set of solutions of any one of the equations
in the system. This approach usually saves considerable amount of
time and effort in trying to solve the system by other methods. In
fact, in some cases this can be the only viable method for solving
the system.

For example, the following system of Diophantine equations (where
$x,y\in\mathbb{Z}$):\\
\begin{tabular*}{16.3cm}{@{\extracolsep{\fill}}lllll}
\noalign{\vskip0.1cm}
$15x+13xy-20y=0$ & and & $3x^{4}+2y^{3}+202=0$ & and & $7x^{5}-4y^{2}-9x^{3}y=484$\tabularnewline[0.1cm]
\end{tabular*}\\
can be easily solved by testing the solutions of the first equation
(assuming we have these solutions or they can be obtained rather easily)
which are:\\
\begin{tabular*}{16.3cm}{@{\extracolsep{\fill}}lll}
\noalign{\vskip0.1cm}
 & $(x,y)=(2,-5)$, $(0,0)$, $(-10,-1)$ & \tabularnewline[0.1cm]
\end{tabular*}\\
on the other equations in the system. By doing so we find out that
only the solution $(x,y)=(2,-5)$ satisfies all the equations of the
system and hence only this solution is acceptable as a solution to
the system.

Similarly, the following system of Diophantine equations:

\begin{tabular*}{15.95cm}{@{\extracolsep{\fill}}llllllll}
\noalign{\vskip0.1cm}
 & $x^{2}-3y=19$ & and & $13y^{3}+6x=11$ & and & $x^{2}+y^{2}=17$ & ($x,y\in\mathbb{Z}$) & \tabularnewline[0.1cm]
\end{tabular*}

\noindent can be easily solved by solving the last equation whose
solutions\footnote{The solutions of the last equation are: $(x,y)=(\pm1,4)$, $(\pm1,-4)$,
$(\pm4,1)$ and $(\pm4,-1)$.} can be easily obtained (because there are only a few possibilities
for $x$ and $y$ to satisfy this equation) and testing these solutions
on the other two equations in the system. On doing so we find that
only the solution $(x,y)=(4,-1)$ satisfies the other two equations
and hence this is the only solution to the system.

\subsubsection{\label{subSubstitution}Substitution, Elimination and Comparison}

Substitution is one of the most common and straightforward methods
for solving systems of Diophantine equations (whether linear or non-linear).
For example, the following system:

\begin{tabular*}{15.95cm}{@{\extracolsep{\fill}}llllll}
\noalign{\vskip0.1cm}
 & $3x^{2}+4y=19$ & and & $5x-2y=3$. & ($x,y\in\mathbb{Z}$) & \tabularnewline[0.1cm]
\end{tabular*}

\noindent can be solved easily by substituting $2y=5x-3$ (which is
obtained from the second equation) into the first equation to get
$3x^{2}+10x-25=0$ which is a univariate quadratic equation in $x$.
Solving this quadratic equation (e.g. by the quadratic formula) will
lead to the only solution of the system which is $(x,y)=(-5,-14)$.

Sometimes this method may require some minor (extra) manipulation
prior to substitution (such as raising to power) to facilitate substitution.
For example, the following system:

\begin{tabular*}{15.95cm}{@{\extracolsep{\fill}}lllll}
\noalign{\vskip0.1cm}
 & $\sqrt{x}-11y^{2}=8$ & and & $4xy-1033y^{2}=411$ & \tabularnewline[0.1cm]
\end{tabular*}

\noindent can be solved by substituting $x=(11y^{2}+8)^{2}$ (which
is obtained from the first equation by squaring $\sqrt{x}$) into
the second equation to get $484y^{5}+704y^{3}-1033y^{2}+256y=411$
whose only integer solution is $y=1$ (i.e. $x=19^{2}=361$) and hence
the only solution to the given system is $(x,y)=(361,1)$.

The methods of elimination and comparison are very similar to the
method of substitution (in their technicalities as well as in their
wide applicability) and hence we do not need to investigate them.

\clearpage{}

\section{Testing the Final Solutions}

It is strongly recommended to test the obtained solutions (whether
of individual equations or of systems of equations) on the given equations
and systems, e.g. by substituting the values of the variables of these
solutions in the given equations and systems. This is particularly
important in the following cases:

\noindent $\bullet$ When the solutions are obtained with certain
manipulations that can introduce foreign solutions (see $\S$ \ref{secManipulationsTransformation}).

\noindent $\bullet$ When the solutions are obtained by messy arguments
and formulations and hence it is likely that the solutions are wrong
or they contain errors and mistakes.

We may also consider using computational tools to do final checks,
e.g. by running a code to obtain the solutions (as we did in $\S$
\ref{secInitialComputInv}) and compare them to the already obtained
solutions. In fact, we strongly recommend using such computational
tools (like coding or spreadsheets or software packages) to check
the final answer especially when we have some doubt or when the produced
argument (or proof or formulation or $\ldots$) is very messy and
susceptible to errors and mistakes. So, we recommend \textit{starting}
our investigation of Diophantine problems by computational inspection
to probe the solution and form a general idea about it (see $\S$
\ref{secInitialComputInv}), and \textit{terminating} our investigation
by computational testing and checking of the obtained solution.

\pagebreak{}

\phantomsection
\addcontentsline{toc}{section}{References}\bibliographystyle{plain}
\bibliography{Bibl}

\clearpage{}

\phantomsection \addcontentsline{toc}{section}{Nomenclature}

\section*{\label{Nomenclature}Nomenclature}

In the following list, we define the common symbols, notations and
abbreviations which are used in the paper as a quick reference for
the reader.\vspace{-0.5cm}\\
\begin{longtable}[l]{ll}
! & factorial\tabularnewline
\noalign{\vskip0.05cm}
$\in$    & in (or belong to)\tabularnewline
\noalign{\vskip0.05cm}
$\ni$ & (backward) in (or belong to)\tabularnewline
\noalign{\vskip0.05cm}
$|a|$ & absolute value of $a$\tabularnewline
\noalign{\vskip0.05cm}
Eq., Eqs. & Equation, Equations\tabularnewline
\noalign{\vskip0.05cm}
LHS, RHS & left hand side, right hand side\tabularnewline
\noalign{\vskip0.05cm}
$m\stackrel{k}{=}n$ & $m$ and $n$ are congruent modulo $k$\tabularnewline
\noalign{\vskip0.05cm}
$m\stackrel{k}{\neq}n$ & $m$ and $n$ are not congruent modulo $k$\tabularnewline
\noalign{\vskip0.05cm}
mod & modulo (or modulus)\tabularnewline
\noalign{\vskip0.05cm}
$\mathbb{N}$ & the set of natural numbers (i.e. $1,2,3,\ldots$)\tabularnewline
\noalign{\vskip0.05cm}
$\mathbb{N}^{0}$ & the set of non-negative integers (i.e. $0,1,2,3,\ldots$)\tabularnewline
\noalign{\vskip0.05cm}
$O_{k}n$ & the order of integer $n$ (modulo $k$)\tabularnewline
\noalign{\vskip0.05cm}
$p$ & prime number\tabularnewline
\noalign{\vskip0.05cm}
$\mathbb{P}$ & the set of prime numbers\tabularnewline
\noalign{\vskip0.05cm}
$\mathbb{Z}$ & the set of integers\tabularnewline
\noalign{\vskip0.05cm}
$\Delta$ & discriminant of quadratic polynomial\tabularnewline
\noalign{\vskip0.05cm}
\end{longtable}

\end{document}